\documentclass[12pt,a4paper]{article}
\usepackage[T1]{fontenc}
\usepackage[utf8]{inputenc}
\usepackage[english]{babel}
\usepackage{lmodern,microtype}
\usepackage[left=28mm,right=28mm,top=38mm,bottom=40mm,headheight=16pt,headsep=9mm]{geometry}
\usepackage{setspace}
\usepackage{amsmath,amssymb,amsthm,mathtools,mathrsfs}
\numberwithin{equation}{section}
\newtheoremstyle{plainspaced}{12pt}{12pt}{\itshape}{}{\bfseries}{.}{.6em}{}
\theoremstyle{plainspaced}
\usepackage{aliascnt}
\newtheorem{theorem}{Theorem}[section]
\newtheorem*{integralrealizationagain}{Theorem~\ref{cc:thm:main}}
\newaliascnt{proposition}{theorem}
\newtheorem{proposition}[proposition]{Proposition}
\aliascntresetthe{proposition}
\newaliascnt{lemma}{theorem}
\newtheorem{lemma}[lemma]{Lemma}
\aliascntresetthe{lemma}
\newaliascnt{corollary}{theorem}
\newtheorem{corollary}[corollary]{Corollary}
\aliascntresetthe{corollary}
\newtheoremstyle{uprightspaced}{12pt}{12pt}{\normalfont}{}{\bfseries}{.}{.6em}{}
\theoremstyle{uprightspaced}

\newtheorem{example}[theorem]{Example}
\newtheorem{question}[theorem]{Question}
\theoremstyle{uprightspaced}
\newaliascnt{remark}{theorem}
\newtheorem{remark}[remark]{Remark}
\newtheorem*{unnumberedremark}{Remark}
\aliascntresetthe{remark}
\usepackage{titlesec}
\titleformat{\section}{\normalfont\large\scshape\centering}{\thesection.}{.65em}{}
\titlespacing*{\section}{0pt}{2.2\baselineskip}{.9\baselineskip}
\titleformat{\subsection}{\normalfont\normalsize\bfseries}{\thesubsection.}{.6em}{}
\titlespacing*{\subsection}{0pt}{1.4\baselineskip}{.5\baselineskip}
\usepackage{fancyhdr}
\newcommand{\shorttitle}{Higher Integral Riemann--Roch}
\usepackage{enumitem}
\setlist[enumerate]{label=\textup{(\roman*)},leftmargin=*,topsep=.5em,itemsep=.3em}
\setlist[itemize]{leftmargin=*,topsep=.5em,itemsep=.3em}
\usepackage[hidelinks,hyperfootnotes=false]{hyperref}
\newcommand{\Z}{\mathbf Z}\newcommand{\Q}{\mathbf Q}\newcommand{\Fp}{\mathbf F_p}
\newcommand{\Zp}{\Z_{(p)}}
\newcommand{\T}{\mathsf T}\newcommand{\one}{\mathbf1}
\newcommand{\SH}{\operatorname{SH}}\newcommand{\CH}{\operatorname{CH}}
\newcommand{\SHeff}{\SH^{\mathrm{eff}}(k)}
\newcommand{\KGL}{\mathrm{KGL}}\newcommand{\kgl}{\mathrm{kgl}}
\newcommand{\HZ}{\mathrm H\Z}\newcommand{\HQ}{\mathrm H\Q}
\newcommand{\cK}{\mathcal K}\newcommand{\Gr}{\operatorname{Gr}}

\newcommand{\Thsp}{\operatorname{Thsp}}
\newcommand{\Th}{\operatorname{Th}}\newcommand{\Td}{\operatorname{Td}}
\newcommand{\ThomClass}[1]{t^{#1}}
\newcommand{\BGL}{\mathrm{BGL}}
\DeclareMathOperator{\cofib}{cofib}\DeclareMathOperator{\fib}{fib}
\DeclareMathOperator{\ch}{ch}\DeclareMathOperator{\rk}{rk}\DeclareMathOperator{\id}{id}
\DeclareMathOperator{\Spec}{Spec}
\DeclareMathOperator{\lcm}{lcm}

\newcommand{\Gm}{\mathbf G_m}
\newcommand{\HF}{\mathrm H\mathbf F}
\usepackage{needspace}
\allowdisplaybreaks[1]
\usepackage{tikz-cd,booktabs,array}
\newcommand{\ct}{\operatorname{CT}}
\newcommand{\Map}{\operatorname{Map}}
\hypersetup{pdftitle={Higher Integral Riemann--Roch},
pdfauthor={Alberto Navarro Garmendia and Jos\'e Navarro Garmendia},
pdfsubject={A proof from integral stable Chern characters},bookmarksnumbered=true}

\DeclareMathOperator{\Mod}{Mod}
\DeclareMathOperator{\im}{im}
\DeclareMathOperator{\Res}{Res}
\newcommand{\iHom}{\underline{\operatorname{Hom}}}

\newcommand{\paperauthors}{A. Navarro Garmendia and J. Navarro Garmendia}
\newcommand{\paperemails}{%
\textsc{Alberto Navarro Garmendia}\\
Departamento de Matem\'atica Aplicada, Universidad Polit\'ecnica de Madrid,\\
Madrid, Spain\\
\href{mailto:alberto.navarro@upm.es}{\texttt{alberto.navarro@upm.es}}\\[1em]
\textsc{Jos\'e Navarro Garmendia}\\
Departamento de Matem\'aticas, Universidad de Extremadura,\\
Badajoz, Spain\\
\href{mailto:navarrogarmendia@unex.es}{\texttt{navarrogarmendia@unex.es}}}

\begin{document}
\begin{center}
{\Large\bfseries Higher Integral Riemann--Roch\par}
\vspace{1em}
{\small\paperauthors\par}
\end{center}
\thispagestyle{plain}
\begin{abstract}
We prove an integral Riemann--Roch formula relating higher algebraic $K$-theory to higher Chow groups for projective morphisms between smooth quasi-projective schemes over a field. The formula is valid after inverting the exponential characteristic. In characteristic zero and nonnegative relative dimension, its specialization to $K_0$ recovers Pappas's formula. 

The proof constructs stable Chern characters on
finite slice truncations of connective algebraic $K$-theory.  

We also present a construction over $\Spec\Z$ that gives an integral formula comparing the oriented Gysin maps in homotopy invariant $K$-theory and motivic cohomology for proper smoothable local complete intersection morphisms of constant virtual dimension.
\end{abstract}



\setcounter{tocdepth}{1}
\begingroup
\small
\setstretch{1}
\setlength{\parskip}{0pt}
\titleformat{\section}{\normalfont\normalsize\scshape\centering}{\thesection.}{.65em}{}
\titlespacing*{\section}{0pt}{.5\baselineskip}{.3\baselineskip}
\makeatletter
\renewcommand*\l@section{\@dottedtocline{1}{0pt}{1.6em}}
\makeatother
\tableofcontents
\endgroup

\section{Introduction}\label{sec:intro}

Grothendieck's insight was to formulate Riemann--Roch as a theorem about
direct images. Let $f\colon X\to Y$ be a projective morphism between smooth
quasi-projective varieties over a field $k$. His theorem says that the square
\[
\begin{tikzcd}[column sep=4.5em,row sep=2.4em]
 K_0(X) \arrow[r,"f_!"] \arrow[d,"{\Td(T_X)\ch}"']
   & K_0(Y) \arrow[d,"{\Td(T_Y)\ch}"] \\
 \CH^*(X)_{\Q} \arrow[r,"f_*"']
   & \CH^*(Y)_{\Q}
\end{tikzcd}
\]
commutes \cite{BorelSerre,SGA6}. Thus for every $x\in K_0(X)$,
\begin{equation}\label{eq:intro-rational-grr}
 \Td(T_Y)\,\ch(f_!x)=f_*\bigl(\Td(T_X)\,\ch(x)\bigr).
\end{equation}
Here $f_!$ is induced by derived direct image of perfect complexes and
$f_*$ is proper pushforward of cycles. The Todd class measures the change
in the Chern character under direct image. Both characteristic series
have rational coefficients. Once their denominators have been cleared in
a given degree, \eqref{eq:intro-rational-grr} tells us that the difference
between the two integral classes is torsion. What becomes of that torsion
is the question of integral Riemann--Roch. We prove such identities for
every higher $K$-group in characteristic zero. The proof determines which
multiples of each Chern character component admit an integral stable lift.

For closed immersions, Grothendieck asked whether the Chern classes of a
direct image could be expressed by an integral polynomial involving the
normal bundle \cite[Exp.~XIV, \S3]{SGA6}. Jouanolou proved this
\emph{Riemann--Roch without denominators} \cite[Th\'eor\`eme~4.1]{JouanolouRR}.
Fulton and MacPherson found integral relations of a different kind for
finite \'etale coverings \cite{FultonMacPherson}. For general projective
morphisms in characteristic zero, Pappas's integral
Riemann--Roch theorem \cite{Pappas} gives a formula in each codimension.
We recall it explicitly, since it is the point of departure for our
formula in higher $K$-theory.

Write $[\,\cdot\,]_j$ for the homogeneous part of codimension $j$. The
Todd number $T_j$ is the exact denominator of the $j$th Todd polynomial
\cite[Lemma~1.7.3]{Hirzebruch}. For a virtual bundle $V$ put
\[
 T_j=\prod_{p\ {\rm prime}}p^{\lfloor j/(p-1)\rfloor},\qquad
 D_j(V)=T_j\Td_j(V),\qquad T_0=D_0(V)=1.
\]
Let $s_i$ be the $i$th Newton polynomial, with $s_0=\rk$. Then
$\ch=\sum_{i\ge0}s_i/i!$ and $\Td=\sum_{j\ge0}D_j/T_j$.
The integral polynomial for their product is Pappas's
\begin{equation}\label{eq:intro-ct}
 \ct_j(x,V)=T_j[\ch(x)\Td(V)]_j
   =\sum_{i=0}^j\frac{T_j}{i!\,T_{j-i}}\,
                   s_i(x)D_{j-i}(V).
\end{equation}
The coefficients on the right are integers \cite[Lemma~2.1 and \S2.b]{Pappas}.
In particular $\ch(x)\Td(V)=\sum_{j\ge0}\ct_j(x,V)/T_j$ as a rational
characteristic series. An expression such as \eqref{eq:intro-ct} always
means that we evaluate the integral polynomial on the right, including
on varieties whose Chow groups have torsion.

Put $d=\dim X-\dim Y$, $T_f=[T_X]-[f^*T_Y]$, and $m=n+d$, where
$n,m\ge0$. Over any field of characteristic zero, the relative forms
of Pappas's identities in $\CH^n(Y)$ are
\begin{equation}\label{eq:intro-pappas}
 \begin{aligned}
  \frac{T_m}{n!}\,s_n(f_!x)
    &=f_*\ct_m(x,T_f) && (d\ge0),\\[2pt]
  \frac{T_n}{n!}\,s_n(f_!x)
    &=\frac{T_n}{T_m}\, f_*\ct_m(x,T_f) && (d<0).
 \end{aligned}
\end{equation}
The first identity is \cite[Corollary~2.3]{Pappas}. For $d<0$,
Pappas states his result with the two absolute tangent bundles, in the form
\(\ct_n(f_!x,T_Y)=(T_n/T_m)f_*\ct_m(x,T_X)\)
\cite[Theorem~2.2(b)]{Pappas}. Summing the identities through degree $n$
after multiplying by the corresponding components of $\Td(-T_Y)$ and using the
projection formula gives the second line of \eqref{eq:intro-pappas}.
The divisibilities $T_iT_{n-i}\mid T_n$ make this an integral
convolution, as detailed in \cite{NavarroGarmendiaIntegral}.

The denominators affect the proof as well as the statement. Factoring
$f$ as a closed immersion followed by a projective bundle projection
introduces an intermediate dimension. Applying the integral formula to
each factor can therefore give a larger multiplier than the one required
for $f$ \cite[Proposition~2.5 and Remark~2.6]{Pappas}. Pappas's proof uses
resolution of singularities and weak factorization. Its final reduction
from a vector bundle to line bundles after blow-ups uses additivity in
$K_0$ \cite[\S6]{Pappas}. For higher $K$-groups there are already
Riemann--Roch theorems for closed immersions
\cite{Gillet,KondoYasuda,NavarroHigherNoDen}. The question we pursue is
how to obtain an identity of the form \eqref{eq:intro-pappas} for a
general projective morphism and a class in any $K_s(X)$.

Our first approach to this question began with the topological
filtration on $K_0$. Only weights through $m$ on $X$ and through $n$
on $Y$ occur in the degree-$n$ formula. This led us in
\cite{NavarroGarmendiaAway} to work with the finite window
\[
 K_{\le r}(X)=K_0(X)/F_{\rm top}^{\,r+1}K_0(X),
\]
where $F_{\rm top}^{\,i}$ is generated by classes of coherent sheaves
supported in codimension at least $i$. Over $\Lambda=\Z[1/r!]$ the
truncated Chern character gives an isomorphism of filtered rings
\[
 (s_i/i!)_{0\le i\le r}\colon K_{\le r}(X)_\Lambda
   \xrightarrow{\ \sim\ }\bigoplus_{i=0}^r\CH^i(X)_\Lambda ,
\]
with the Chow product truncated above $r$. Comparing direct images on
the source and target windows then gives the degree-$n$ Riemann--Roch
formula over $\Z[1/N!]$, with $N=\max\{n,m+1\}$, over an arbitrary
field. In the companion paper \cite{NavarroGarmendiaIntegral}, we
retain the integral filtration. After clearing the two sides with
$D=\lcm(n!,T_m)$, algebraic Adams--Riemann--Roch gives an explicit
annihilator $V_m$ for the resulting defect. On $F_{\rm top}^{\,i}K_0(X)$
the bound improves to $V_{m-i}$ for $0\le i\le m$.

The idea of the present paper is to pass from a finite range of
classes to a finite range of the spectrum that represents them.
In this form a single integral operation can see every higher
$K$-group. Its value on a Thom class also records what happens
under direct image. We first state the formula and then describe
the operation from which it follows.

\medskip
\noindent\emph{A note on this version.}
We used ChatGPT and Claude to adapt the calculation of Dundas and
\O{}stv{\ae}r in Section~\ref{sec:integral-realization}. We have
checked the resulting computations and believe they are correct, but
in normal circumstances we would have waited to discuss them further
with colleagues. In Section~\ref{sec:arithmetic-rr} we give only an
outline of the argument, adapting the ideas of
Section~\ref{sec:integral-realization} to the arithmetic setting
using the work of Bachmann and Hopkins. Following advice from
colleagues about the confidentiality of work discussed with commercial
AI systems, we have chosen to circulate this preprint earlier than we
ordinarily would.

\subsection{The Riemann--Roch formula}

Let $X$ and $Y$ be smooth quasi-projective $k$-schemes of pure
dimension and $f\colon X\to Y$ projective. Keep $d$, $T_f$ and $m=n+d$
as above. Write $f_!\colon K_s(X)\to K_s(Y)$ for the direct image
induced by $Rf_*$ \cite[\S3.16.4]{TT}, and
$f_*\colon\CH^j(X,s)\to\CH^{j-d}(Y,s)$ for proper pushforward in
Bloch's higher Chow groups \cite{Bloch}. The Newton operations are
\begin{equation}\label{eq:newton}
 \chi_{j,0}=s_j,\qquad
 \chi_{j,s}=(-1)^{j-1}j\,c_{j,s}\quad(j>0,\ s>0),\qquad
 \chi_{0,s}=0\quad(s>0).
\end{equation}
The second formula follows from Newton's identity. On a simplicial
suspension, products of two reduced classes vanish, so only the
linear term of $s_j$ survives. For a characteristic series $Q$ independent
of $x$, set
\begin{equation}\label{eq:linear}
 [\ch(x)Q]^{(s)}_j=\sum_{i=0}^j\frac{\chi_{i,s}(x)}{i!}Q_{j-i}.
\end{equation}

\Needspace{11\baselineskip}
\begin{theorem}[Higher integral Riemann--Roch]\label{thm:main}
Suppose $k$ has characteristic zero. Let $n,s\ge0$ and $x\in K_s(X)$.
If $m=n+d\ge0$, put
\begin{equation}\label{eq:multiplier}
 C_{n,d}=\begin{cases}
 \lcm(n!,T_m),& d\ge-2\ \text{or}\ m=0,\\
 \lcm(n!,T_{n-2}),& d\le-3,\ m>0 .
 \end{cases}
\end{equation}
Then in $\CH^n(Y,s)$ one has
\begin{equation}\label{eq:main}
 \frac{C_{n,d}}{n!}\,\chi_{n,s}(f_!x)
   =f_*\!\left(C_{n,d}[\ch(x)\Td(T_f)]^{(s)}_m\right).
\end{equation}
If $m<0$, then $\chi_{n,s}(f_!x)=0$.
\end{theorem}

Both sides are evaluated integrally. In particular, the class inside
the pushforward is
\begin{equation}\label{eq:integral-expression}
 \sum_{i=0}^m\frac{C_{n,d}}{i!\,T_{m-i}}\,
            \chi_{i,s}(x)D_{m-i}(T_f).
\end{equation}
For $s=0$ and $d\ge0$ we recover the first line of
\eqref{eq:intro-pappas}, with $C_{n,d}=T_m$. For $d<0$ our
multiplier divides $T_n$ and can be smaller. For example,
$C_{3,-2}=6$, whereas $T_3=24$. In source weight zero,
$C_{n,-n}=n!$. For $s>0$ and $n>0$ the formula in higher Chern
classes reads, with $C=C_{n,d}$,
\begin{equation}\label{eq:higher-explicit}
 (-1)^{n-1}\frac{C}{(n-1)!}c_{n,s}(f_!x)
 = f_*\!\left(\sum_{i=1}^m(-1)^{i-1}\frac{C}{(i-1)!}
                    \Td_{m-i}(T_f)c_{i,s}(x)\right).
\end{equation}
The empty sum is zero. For $s>0$ and $n=0$, $\CH^0(Y,s)=0$.

Already in codimension one, the factor on the left has geometric
content. Take an elliptic curve $Y$ over an algebraically closed field
of characteristic zero and a nontrivial line bundle $L$ of order two.
A trivialization of $L^{\otimes2}$ gives a finite \'etale double cover
$f\colon X=\Spec_Y(\mathcal O_Y\oplus L)\to Y$. Then $T_f=0$ and
$c_1(f_!\mathcal O_X)=c_1(L)\ne0$, while $f_*c_1(\mathcal O_X)=0$.
Thus the factor $T_1=2$ is necessary. In higher weights, the optimal
factor for a particular class of direct images can be smaller than the
one supplied by a stable operation. For finite \'etale covers in
characteristic zero, the optimal universal multiplier on the $n$th
Newton class of a $K_0$-class is the product
of the distinct primes dividing $T_n/n!$
\cite[Theorem~23.3]{FultonMacPherson}, \cite[\S2.d.1]{Pappas}.

\subsection{Stable integrality and the proof}

The spectral analogue of the filtration used on $K_0$ is the slice
filtration of algebraic $K$-theory \cite{Voevodsky,Levine}. Write
$\T=\Sigma^{2,1}$ and let $\kgl=f_0\KGL$ be the effective cover of
periodic algebraic $K$-theory. Its window of weights $0,\ldots,q$ and
its slice of weight $q$ are, respectively,
\[
 W_q=\cofib\bigl(\T^{q+1}\kgl
            \xrightarrow{\ \beta^{q+1}\ }\kgl\bigr),
 \qquad B_q=\T^q\HZ .
\]
Here $\beta$ is the Bott map and we continue to work over a field of
characteristic zero. Set $H_M=\bigoplus_{q=0}^M B_q$. Rationally the
window splits into its weights, and we denote the projections by
$\ch_q$ with the normalization of Riou
\cite[Definition~6.2.3.9]{Riou}. For $q\ge1$ the integral attaching map
$\kappa_q\colon W_{q-1}\to\Sigma B_q$ measures the obstruction to
extending the identity of the top slice across $W_q$.

\begin{theorem}[Stable integral Chern character]\label{thm:intro-stable}
For $q\ge0$, the image of $[\kgl,B_q]$ in
$[\kgl,\T^q\HQ]=\Q\cdot\ch_q$ is $\Z\cdot T_q\ch_q$. There is
a map $\theta_q\colon W_q\to B_q$ whose restriction to $B_q$ is
$T_q\id$ and whose rationalization is $T_q\ch_q$. For $q\ge1$,
$\kappa_q$ has order $T_q$. Let $\pi_{M,q}\colon W_M\to W_q$ be
the quotient. The components
\[
 (\Phi_M)_q=\frac{T_M}{T_q}\,\theta_q\circ\pi_{M,q}
 \qquad(0\le q\le M)
\]
define a map $\Phi_M\colon W_M\to H_M$ with rationalization
$T_M\ch_{\le M}$. It acts by $T_M$ on each slice and becomes an
equivalence after inverting $T_M$. In weight $q$ it induces
$(T_M/q!)\chi_{q,s}$ on $K_s(X)$, independently of the choices of
the $\theta_q$.
\end{theorem}

The same Todd numbers occur in the stable Chern characters of
connective complex $K$-theory studied by Adams \cite{Adams,Anderson}.
To construct $\theta_q$ motivically we work one prime at a time.
Dundas and \O{}stv{\ae}r compute the integral motivic dual Steenrod
algebra. In particular they describe the Bockstein on cooperations
and show that $p$ kills the reduced $p$-local cooperations
\cite[Remark~1.1, Lemma~1.2 and Theorem~1.4]{DO}. Hoyois describes
$\HF_p\wedge\HZ$ as a sum of suspensions of $\HF_p$ as an
$\HF_p$-module \cite[Theorem~5.17]{Hoyois}. The step we need is to
realize the pairs in the Bockstein calculation by maps of
$\HZ_{(p)}$-modules. This gives the following description.

\begin{theorem}[Integral motivic cooperations]\label{thm:intro-cooperations}
For every prime $p$ there is an equivalence of left $\HZ_{(p)}$-modules
\[
 \HZ_{(p)}\wedge\HZ_{(p)}
 \simeq \HZ_{(p)}\oplus\bigoplus_\lambda
    \Sigma^{2w_\lambda+\ell_\lambda-1,w_\lambda}\HF_p,
 \qquad w_\lambda\ge p-1,\quad \ell_\lambda\ge1 .
\]
In particular, $p$ annihilates
$[\HZ_{(p)},\Sigma^{2h+1,h}\HZ_{(p)}]$ for every $h\ge0$.
This group vanishes for $h<p-1$.
\end{theorem}

These odd operations are precisely where the extension obstructions
take their values. Adams operations separate weights modulo $p-1$.
In the summand containing weight $q$, the lower weights are
$q-(p-1),q-2(p-1),\ldots$ down to the last nonnegative term.
At each step a factor of $p$ kills the obstruction. Gluing the local
extensions gives $\theta_q$ with top value
$\prod_p p^{\lfloor q/(p-1)\rfloor}=T_q$. Evaluation on universal
Thom classes shows that every stable lift must have a top value
divisible by this same integer (Theorem~\ref{thm:exact-denominator}).

To see how the stable character yields Riemann--Roch, consider the
proper transfer
\[
 \tau_f\colon\Sigma^\infty Y_+\longrightarrow\Th(-T_f).
\]
Duality gives a description of both direct images by precomposition
with this one map, after the appropriate Thom identifications
(Lemma~\ref{lem:transfer-comparison}). This viewpoint underlies
\cite{NavarroNavarro} and is developed for smooth proper schemes in
\cite{NavarroGarmendiaDuality}. A stable operation $\psi$ acts by
postcomposition. For a class $z$ on the Thom spectrum, associativity
gives $(\psi\circ z)\circ\tau_f=\psi\circ(z\circ\tau_f)$, with the
corresponding suspension for a class of higher degree. To read this
equality in Chow groups one must compute $\psi$ on the product of
a $K$-class with the Thom class of $-T_f$. Rationally the Thom
correction is $\Td(-T_f)^{-1}=\Td(T_f)$. The relevant integral
cohomology of universal Grassmannians and their Thom spaces is
torsion free, so this rational computation determines the integral
polynomial there, before it is pulled back to $X$. This is the
calculation in Proposition~\ref{prop:thom}. It applies to the
appropriate Bott translates of $\Phi_M$ and gives the formula
with $T_M$ for $M=\max\{n,m\}$. For $d<0$, a shorter window and
a separate argument in source weight zero give \eqref{eq:multiplier}.

Sections~\ref{sec:preliminaries} and~\ref{sec:integral-realization}
develop the stable construction, and Section~\ref{sec:transfer} gives
the transfer argument in characteristic zero.
Section~\ref{sec:characteristic} treats Bloch's higher Chow groups
over a perfect field of exponential characteristic $e$ with
$\Z[1/e]$-coefficients, and then descends the result to an arbitrary
field. The multipliers are those of Theorem~\ref{thm:main}.
Section~\ref{sec:arithmetic-rr} constructs the character over
$\Spec\Z$ itself. The calculations of Bachmann and Hopkins include
the residue characteristic \cite[Theorems~7.1--7.2]{BachmannHopkinsCharacteristic}.
The base-change results of Bachmann, Elmanto and Morrow
\cite[Corollaries~9.6--9.7]{BEM} allow us to pull back the construction.
For a proper morphism of quasi-compact quasi-separated schemes admitting
a factorization as a regular closed immersion of finite presentation
followed by a smooth separated morphism of finite presentation, and
with $T_f$ of constant rank, Theorem~\ref{thm:arithmetic-rr} gives an
integral identity from $KH_s(X)$ to $H^{2n-s}(Y,\Z(n))$ with the
oriented Gysin maps and multiplier $T_{\max(n,m)}$.
Corollary~\ref{cor:arithmetic-smooth} recovers the multipliers
and boundary cases of Theorem~\ref{thm:main} for projective morphisms
between smooth quasi-projective $\Z$-schemes of pure relative dimensions. The comparison used here with proper pushforward on Bloch's
higher cycles over a field inverts $e$
\cite[Corollary~8.12]{CDIntegral}.
Section~\ref{sec:examples} studies the coefficients for individual
direct images.

\section{Preliminaries}\label{sec:preliminaries}

For a virtual bundle $\xi$, put $D_j(\xi)=T_j\Td_j(\xi)$ as in
\cite[\S\S2.b--2.c]{Pappas}. We recall Hirzebruch's formula for the
denominators of the Todd polynomials \cite[Lemma~1.7.3]{Hirzebruch}
and the divisibilities given by \cite[Lemma~2.1]{Pappas}.

\begin{lemma}\label{lem:denominators}
The $D_j$ are integral Chern polynomials. Moreover,
\begin{equation}\label{eq:divisibility}
 i!T_{j-i}\mid T_{j-1}\quad(1\le i\le j),\qquad
 T_aT_b\mid T_{a+b},\qquad (j+1)!\mid T_j.
\end{equation}
In particular, $T_j\mid T_M$ for $j\le M$ and
$\Z[1/T_j]=\Z[1/(j+1)!]$.
\end{lemma}

Consequently, for $M\ge j$ the expression
$T_M[\ch(x)\Td(\xi)]^{(s)}_j$ is the integral operation
\begin{equation}\label{eq:integral-polynomial}
 \sum_{i=0}^j\frac{T_M}{i!T_{j-i}}\chi_{i,s}(x)D_{j-i}(\xi).
\end{equation}
For $s=0$ and $M=j$, this is Pappas's polynomial $\ct_j$
\cite[\S2.b]{Pappas}. The same convention applies with $\Td(\xi)^{-1}=\Td(-\xi)$.
For the second branch of \eqref{eq:multiplier}, $n-2\ge m$;
in both branches $C_{n,d}$ is divisible by $n!$ and $T_m$. The case $m=0$ is immediate. Thus all
expressions in the theorem are integral.

\subsection{Conventions}
Let $p\colon X\to\Spec k$ be smooth. For a vector bundle $V$ on $X$,
write $\Th_X(V)\in\SH(X)$ for its Thom object and
$\Th(V)=p_{\#}(\Th_X(V))\in\SH(k)$. For a virtual bundle $\xi=E-F$,
we use $\Th_X(\xi)\simeq\Th_X(E)\otimes\Th_X(F)^{-1}$
\cite[\S\S2.1.4--2.1.5]{DJK}. Write $\T=\Sigma^{2,1}$,
$\Sigma=\Sigma^{1,0}$ and $S^s=(S_s^1)^{\wedge s}$. All smash products
are over the motivic sphere unless otherwise indicated.
We use the standard orientations on the commutative ring spectra $\KGL$
and $\HZ$. For periodic algebraic $K$-theory the first Chern class of a
line bundle $L$ is $c_1^{\KGL}(L)=\beta^{-1}(1-[L^\vee])$
\cite[Examples~2.1.4(4) and~2.1.9]{DegliseOrientation}.
For integral motivic cohomology we use Voevodsky's motivic
Eilenberg--Mac Lane spectrum $\HZ$ over $k$, constructed using finite
correspondences, with the additive orientation given by the usual first
Chern class \cite[\S\S4.1--4.3, in particular~(4.27)]{Hoyois}.
Let $\HQ=(\HZ)_\Q$ denote the spectrum obtained by extension of
coefficients to $\Q$. Its multiplication
$\HQ\wedge\HQ\to\HQ$ is an equivalence
\cite[Proposition~14.1.6 and Corollary~16.1.7(1)]{CDMotives}.
We write $[E,F]=\operatorname{Hom}_{\SH(k)}(E,F)$ and denote by
$\Map(E,F)$ the mapping spectrum in the underlying stable
$\infty$-category, so that $[E,F]=\pi_0\Map(E,F)$.

For the slice filtration, let $\SHeff$ be the localizing subcategory of
$\SH(k)$ generated by $\Sigma^\infty U_+$ for $U$ smooth over $k$.
The $i$th effective cover $f_iE\to E$ is universal among maps to $E$
from objects of $\T^i\SHeff$, and the $i$th slice is
$s_iE=\cofib(f_{i+1}E\to f_iE)$ \cite[\S8.3]{Hoyois}.
An endomorphism of $E$ thus acts compatibly on its covers, slices and
finite truncations $\cofib(f_bE\to f_aE)$ for $b>a$.
The subcategories $\T^i\SHeff$ are stable under all simplicial shifts and
are generated by compact objects of $\SH(k)$, so $f_i$ and $s_i$
preserve colimits. The slice $s_iE$ lies in $\T^i\SHeff$ and satisfies
$[Y,s_iE]=0$ for $Y\in\T^{i+1}\SHeff$.
For $\KGL$ we put
\[
 A_i=f_i\KGL,\qquad B_i=s_i\KGL,\qquad \kgl=A_0.
\]
The inclusion $\SHeff\hookrightarrow\SH(k)$ is symmetric monoidal, so its
right adjoint is lax symmetric monoidal and $\kgl=f_0\KGL$ inherits a
ring structure. More generally, $A_a\wedge A_b$ is $(a+b)$-effective,
so the product $A_a\wedge A_b\to\KGL$ factors canonically through
$A_{a+b}$. We equip $\kgl$ with the standard orientation lifting that
of $\KGL$ along the cover map.
For $X$ smooth over $k$ and $s\ge0$, the source
$\Sigma^{s,0}\Sigma^\infty X_+$ is effective. The cover map therefore
gives the first isomorphism in
\[
 [\Sigma^{s,0}\Sigma^\infty X_+,\kgl]
 \xrightarrow{\ \sim\ }[\Sigma^{s,0}\Sigma^\infty X_+,\KGL]
 \cong K_s(X),
\]
where the second isomorphism is representability of Quillen $K$-theory
on regular schemes \cite[Remark~2.1.5]{DegliseOrientation}.

For a localization $\Delta\subseteq\Q$ of $\Z$, put $E_{\Delta}=E\wedge\one_{\Delta}$, where
$\one_{\Delta}$ is the localization of the sphere obtained by inverting the
integers that are units in $\Delta$. As $E_{\Delta}$ is a sequential colimit of
copies of $E$, one has $f_i(E_{\Delta})\simeq(f_iE)_{\Delta}$ and
$s_i(E_{\Delta})\simeq(s_iE)_{\Delta}$, so localization preserves all finite slice
towers below.

We identify $s_0\one$ with $\HZ$ via the unit of motivic cohomology
\cite[Theorem~10.5.1]{Levine}. The unit $\eta\colon\one\to\KGL$ then
induces the equivalence
$\alpha_0\colon\HZ\simeq s_0\one\xrightarrow{s_0(\eta)}B_0$
\cite[Theorem~4.1]{RondigsOstvaer}.
For $i\in\Z$, multiplication by the $i$th power of the Bott element gives
an equivalence $\beta^i\colon\T^i\KGL\xrightarrow{\sim}\KGL$, where
$\beta$ is represented on $\mathbf P^1$ by $[\mathcal O(1)]-1$.
Since $f_i(\T^i E)\simeq\T^i f_0E$ and $s_i(\T^i E)\simeq\T^i s_0E$,
applying $f_i$ and $s_i$ to $\beta^i$ gives the equivalences
\[
 \T^i\kgl\xrightarrow{\sim}A_i,\qquad
 \alpha_i\colon\T^i\HZ\xrightarrow{\T^i\alpha_0}
 \T^i B_0\xrightarrow{s_i(\beta^i)}B_i.
\]
Write $\iota_{b,a}\colon A_b\to A_a$ for the tower map when $b\ge a$.
Under these identifications $\iota_{i+1,i}$ is the $i$th Tate suspension
of the connective Bott map $\beta\colon\T\kgl\to\kgl$, and its
cofiber is $B_i$.

For each integer $u\ge1$, choose a stable Adams endomorphism $\Psi^u$ of
$\KGL[1/u]$ inducing the usual Adams operation on algebraic $K$-theory,
with $\Psi^1=\id$ \cite[Definition~5.3.2 and the paragraph following it]{Riou}.
We use these endomorphisms in the homotopy category. All the evaluations
below are independent of this choice by Lemma~\ref{lem:evaluation}.
On the rational Adams summand of weight $i$, the operation $\Psi^u$
acts by $u^i$ \cite[Proposition~5.3.14]{Riou}.

The rational Chern character is the equivalence
\[
 \ch=(\ch_i)_{i\in\Z}\colon\KGL_\Q\xrightarrow{\sim}
     \bigoplus_{i\in\Z}\T^i\HQ
\]
whose weight-zero component induces the rank on $K_0$ and whose other
components are obtained from $\ch_0$ by Bott periodicity. It identifies
the rational Adams summand of weight $i$ with $\T^i\HQ$
\cite[Definition~6.2.3.9 and Remark~6.2.3.10]{Riou}.

\needspace{6\baselineskip}
\begin{lemma}\label{lem:slice-normalization}
For $i\in\Z$, the composite
\[
 \T^i\HQ\xrightarrow{\alpha_{i,\Q}} B_{i,\Q}
 \xrightarrow{s_i(\ch_i)}\T^i\HQ
\]
is the identity, where we use $s_i(\KGL_\Q)\simeq B_{i,\Q}$ and
$s_i(\T^i\HQ)\simeq\T^i\HQ$.
\end{lemma}
\begin{proof}
Since the unit has rank one, $\ch_0\circ\eta_\Q$ is the unit of $\HQ$.
Taking its zero slice gives $s_0(\ch_0)\circ\alpha_{0,\Q}=\id_{\HQ}$.
By the definition of the Bott translates,
$\ch_i\circ\beta^i_\Q=\T^i\ch_0$
\cite[Definition~6.2.3.9]{Riou}. In weight one this convention sends
$[\mathcal O(1)]-1$ on $\mathbf P^1$ to $h=c_1(\mathcal O(1))$,
since $\exp(h)-1=h$ there. Applying $s_i$ to the Bott identity and
composing with $\T^i\alpha_{0,\Q}$ gives
$s_i(\ch_i)\circ\alpha_{i,\Q}=\T^i(s_0(\ch_0)\circ\alpha_{0,\Q})$,
as claimed.
\end{proof}

\begin{unnumberedremark}
Let $Z\subset X$ be an integral closed subvariety of codimension $d>0$
in a smooth quasi-projective variety, and put $x=[\mathcal O_Z]\in K_0(X)$.
Riemann--Roch without denominators gives $c_j(x)=0$ for $0<j<d$ and
$c_d(x)=(-1)^{d-1}(d-1)![Z]$
\cite[Theorem~15.3 and Example~15.3.6]{Fulton}. The vanishing of the
lower Chern classes reduces Newton's identity to
\[
 \chi_{d,0}(x)=(-1)^{d-1}d\,c_d(x)=d![Z]\quad\text{in }\CH^d(X).
\]
Thus $\ch_d(x)=\chi_{d,0}(x)/d!=[Z]$ in $\CH^d(X)_{\Q}$.
Bott desuspension sends $\chi_{d,0}$ to $d\,\chi_{d-1,0}$ and hence
$\ch_d$ to $\ch_{d-1}$ \cite[Lemma~6.2.3.2]{Riou}, which explains
the scalar $1$ in Lemma~\ref{lem:slice-normalization}.
\end{unnumberedremark}

\section{Integral realization on finite truncations}\label{sec:integral-realization}
For $M\geq0$ define the finite slice truncation and its split counterpart by
\begin{equation}\label{cc:eq:truncations}
 W_M=\cofib(A_{M+1}\longrightarrow A_0),\qquad
 H_M=\bigoplus_{q=0}^M B_q.
\end{equation}
Thus $W_M\simeq\kgl/\beta^{M+1}$, where $\beta\colon\T\kgl\to\kgl$ is
the Bott map; the quotient notation refers to a cofiber. The slices of $W_M$ are $B_0,\ldots,B_M$.
For $q\geq1$, the top slice of $W_q$ gives a cofiber sequence
\begin{equation}\label{cc:eq:attaching}
 B_q\xrightarrow{\ \iota_q\ }W_q\longrightarrow W_{q-1}
 \xrightarrow{\ \kappa_q\ }\Sigma B_q,
\end{equation}
and we call $\kappa_q$ the \emph{top attaching map}. Set $W_{-1}=0$,
$W_0=B_0=\HZ$ and $\iota_0=\id$. 

\begin{unnumberedremark}
For $q\ge1$, extending $\id_{B_q}$ along $\iota_q$ would give a
retraction $W_q\to B_q$ and split \eqref{cc:eq:attaching} as
$W_q\simeq B_q\oplus W_{q-1}$. More generally, the connecting
homomorphism for $[-,B_q]$ sends $N\id_{B_q}$ to $N\kappa_q$ in
$[W_{q-1},\Sigma B_q]$. Thus $N\id_{B_q}$ extends precisely when
$N\kappa_q=0$.
\end{unnumberedremark}

\newcommand{\IntegralRealizationStatement}[2]{%
Let $q\geq0$.
\begin{enumerate}[label=\textup{(\alph*)}]
\item There is a map $\theta_q\colon W_q\to B_q$ with
 $\theta_q\iota_q=T_q\id_{B_q}$, and one may take $\theta_0=\id_{\HZ}$.
 Equivalently, $T_q\kappa_q=0$ for $q\geq1$.
\item Every map $\theta\colon W_q\to B_q$ with $\theta\iota_q=N\id_{B_q}$
 rationalizes to $N\ch_q$, where $\ch_q\colon W_{q,\Q}\to B_{q,\Q}$ is
 the weight-$q$ Chern projection characterized in Lemma~\ref{cc:lem:rational}. Two choices of $\theta_q$ differ by a
 map that factors through $W_q\to W_{q-1}$ and vanishes after
 rationalization.
\item Fix choices of the $\theta_q$. For $M\geq0$, let
 $H_M=\bigoplus_{i=0}^M B_i$ be the direct sum of the slices of $W_M$.
 The maps
 \begin{#1}#2
  \Phi_M\colon W_M\longrightarrow H_M,\qquad
  (\Phi_M)_q=\frac{T_M}{T_q}\,\theta_q\circ\pi_{M,q},
 \end{#1}
 where $\pi_{M,q}\colon W_M\to W_q$ is the quotient map, satisfy
 \[
  s_q(\Phi_M)=T_M\id_{B_q}\quad(0\leq q\leq M),
  \qquad (\Phi_M)_{\Q}=T_M\ch_{\leq M}.
 \]
 Here $\ch_{\leq M}\colon W_{M,\Q}\to H_{M,\Q}$ is the truncation of
 Riou's Chern character with $q$th component $\ch_q\circ(\pi_{M,q})_{\Q}$.
 The map $\Phi_M$ becomes an equivalence after inverting $T_M$.
\end{enumerate}
}

\begin{theorem}[Integral realization]\label{cc:thm:main}
\IntegralRealizationStatement{equation}{\label{cc:eq:Phi-intro}}
\end{theorem}

The existence proof is given in Subsections~\ref{cc:sec:operations}--\ref{cc:sec:glue}.
Theorem~\ref{thm:exact-denominator} shows that $T_q$ is the least positive
integer $N$ for which $N\ch_q$ lifts to an integral stable operation.

\subsection{Rational rigidity}
Since $f_iA_0=A_i$ for $i\geq0$ and $f_iA_{q+1}=A_{q+1}$ for
$i\leq q+1$, we have
\[
 f_iW_q=\cofib(A_{q+1}\longrightarrow A_i)\qquad(0\leq i\leq q+1).
\]
Thus $f_0W_q=W_q$, $f_{q+1}W_q=0$, and there are cofiber sequences
\begin{equation}\label{cc:eq:filtration}
 f_{i+1}W_q\longrightarrow f_iW_q\longrightarrow B_i
 \longrightarrow\Sigma f_{i+1}W_q\qquad(0\leq i\leq q).
\end{equation}
Moreover $f_qW_q=B_q$, the inclusion $f_qW_q\to W_q$ is $\iota_q$, and
$W_q/f_qW_q\simeq W_{q-1}$; this gives \eqref{cc:eq:attaching}.

\begin{lemma}\label{cc:lem:descent}
Let $M\geq q\geq0$. Restriction along $A_q\to B_q$ and along
$A_0\to W_M$ gives equivalences
\begin{equation}\label{cc:eq:descent}
 \Map(B_q,B_q)\simeq\Map(A_q,B_q),\qquad
 \Map(W_M,B_q)\simeq\Map(A_0,B_q).
\end{equation}
Moreover $[B_q,B_q]=\Z$, and for $q\geq1$ and $N\in\Z$,
\begin{equation}\label{cc:eq:criterion}
 \exists\,\theta\colon W_q\to B_q\text{ with }\theta\iota_q=N\id_{B_q}
 \quad\Longleftrightarrow\quad N\kappa_q=0.
\end{equation}
The same holds after localization at a subring $\Delta\subseteq\Q$, with
$[B_{q,\Delta},B_{q,\Delta}]=\Delta$.
\end{lemma}

\begin{proof}
The fibers of $A_q\to B_q$ and $A_0\to W_M$ are $A_{q+1}$ and $A_{M+1}$.
Both lie in $\T^{q+1}\SHeff$ because $M\geq q$. For either fiber $Y$,
slice orthogonality gives
$\pi_n\Map(Y,B_q)=[\Sigma^{n,0}Y,B_q]=0$ for every $n\in\Z$,
since $\Sigma^{n,0}Y$ is still $(q+1)$-effective.
This proves \eqref{cc:eq:descent}. Since $s_0\one\simeq\HZ$,
the cofiber sequence $f_1\one\to\one\to\HZ$ and slice orthogonality give
$[\HZ,\HZ]\cong[\one,\HZ]=H^0(k,\Z(0))=\Z$, and $\T$ is invertible.
Finally, $[-,B_q]$ applied to \eqref{cc:eq:attaching} gives the exact
sequence
\[
 [W_q,B_q]\xrightarrow{\ \iota_q^*\ }[B_q,B_q]
 \xrightarrow{\ \partial\ }[\Sigma^{-1}W_{q-1},B_q],
 \qquad\partial(\phi)=\phi\circ\Sigma^{-1}\kappa_q,
\]
so $N\id$ lies in the image of $\iota_q^*$ if and only if
$N\Sigma^{-1}\kappa_q=0$. The localized statements have the same proof,
using $s_0(\one_{\Delta})\simeq\HZ_{\Delta}$ and $[\one,\HZ_{\Delta}]=\Delta$.
\end{proof}

Taking $M=q$ in \eqref{cc:eq:descent} gives $[W_q,B_q]\cong[\kgl,\T^q\HZ]$.
Under this identification, a map $\theta$ as in \eqref{cc:eq:criterion}
corresponds to a stable operation whose restriction to $f_q\kgl=A_q$ is
$N$ times the projection $A_q\to B_q$. In particular, the map $\theta_q$
of Theorem~\ref{cc:thm:main}(a) gives an integral stable operation
$\kgl\to\T^q\HZ$ whose rationalization is $T_q\ch_q$ by part~(b).
Operations on $\kgl$ therefore descend uniquely to truncations, and the
substance of Theorem~\ref{cc:thm:main} is the extension problem
\eqref{cc:eq:criterion}.

\needspace{8\baselineskip}
Over $\Q$, the weights below $q$ contribute neither an obstruction to
extending an endomorphism of the top slice nor an ambiguity in its extension.

\begin{lemma}\label{cc:lem:rational}
Let $q\geq0$. Then
\[ [W_{q-1,\Q},\Sigma^{n,0}B_{q,\Q}]=0\qquad(n=0,1), \]
and restriction
along $\iota_q$ is an isomorphism
\[
 [W_{q,\Q},B_{q,\Q}]\xrightarrow{\ \sim\ }[B_{q,\Q},B_{q,\Q}]=\Q,
\]
whose inverse sends $1$ to $\ch_q$.
\end{lemma}

\begin{proof}
For $0\leq i\leq q$, put $j=q-i$ and $D=\Sigma^{2j+n,j}\HQ$.
Cancelling $\T^i$ and using the free $\HQ$-module adjunction and the
multiplication equivalence $\HQ\wedge\HQ\simeq\HQ$, we obtain
\[
 \begin{aligned}
 [B_{i,\Q},\Sigma^{n,0}B_{q,\Q}]
 &\cong [\HQ,D]
  \cong [\HQ\wedge\HQ,D]_{\Mod_{\HQ}}\\
 &\cong [\HQ,D]_{\Mod_{\HQ}}
  \cong [\one,D]
  = H^{2j+n}(k,\Q(j)).
 \end{aligned}
\]
For $n\in\{0,1\}$ this vanishes unless $j=n=0$, because
$H^m(k,\Q(j))=0$ for $m>j$ and for $j<0$. This is one
place where the field hypothesis enters. Induction along the
filtration \eqref{cc:eq:filtration} of $W_{q-1}$, whose slices are
$B_0,\ldots,B_{q-1}$, gives the first claim. The exact sequence
\[
 [W_{q-1,\Q},B_{q,\Q}]\to[W_{q,\Q},B_{q,\Q}]\to[B_{q,\Q},B_{q,\Q}]
 \to[W_{q-1,\Q},\Sigma B_{q,\Q}]
\]
obtained from \eqref{cc:eq:attaching} then gives the isomorphism. The
weight-$q$ component of $\ch$ kills $A_{q+1,\Q}$, so it factors through
$W_{q,\Q}$, and its top value is $1$ by Lemma~\ref{lem:slice-normalization}.
\end{proof}

\subsection{From cooperations to odd stable operations}\label{cc:sec:operations}

Fix a prime $p$ and write
\[
 H=\pi_{**}(\HF_p)=\bigoplus_{a,b\in\Z}H^{-a}(k,\Fp(-b))
\]
for the bigraded ring of motivic cohomology of $k$ with $\Fp$-coefficients,
indexed homologically by $\pi_{a,b}E=[\Sigma^{a,b}\one,E]$.
Bocksteins have bidegree $(-1,0)$ on homotopy groups, since
$\pi_{a,b}(\Sigma^{1,0}E)=\pi_{a-1,b}E$. We
work in the stable $\infty$-category $\Mod_{\HZ_{(p)}}$ of left $\HZ_{(p)}$-modules; a
free module $\HZ_{(p)}\wedge E$ carries the action on the left factor.

Let $u\colon\one\to \HZ_{(p)}$ denote the unit of $\HZ_{(p)}$. The multiplication
$\mu\colon \HZ_{(p)}\wedge \HZ_{(p)}\to \HZ_{(p)}$ is $\HZ_{(p)}$-linear and has $\id_{\HZ_{(p)}}\wedge u$ as
an $\HZ_{(p)}$-linear section. Thus
\begin{equation}\label{cc:eq:C}
 \HZ_{(p)}\wedge \HZ_{(p)}\simeq \HZ_{(p)}\oplus C,\qquad C=\fib(\mu).
\end{equation}
Put $C/p=\HF_p\otimes_{\HZ_{(p)}} C$. It is a left $\HF_p$-module, and there is a cofiber
sequence of $\HZ_{(p)}$-modules
\begin{equation}\label{cc:eq:bockstein-triangle}
 C\xrightarrow{\ p\ }C\xrightarrow{\ r\ }C/p\xrightarrow{\ \delta\ }\Sigma C.
\end{equation}
Since $\HF_p\otimes_{\HZ_{(p)}}(\HZ_{(p)}\wedge \HZ_{(p)})\simeq \HF_p\wedge \HZ_{(p)}$, \eqref{cc:eq:C} gives
$\HF_p\wedge \HZ_{(p)}\simeq \HF_p\oplus C/p$. The Bockstein of $C/p$ is the map
$d=r\circ\delta$ on homotopy, of degree $(-1,0)$.

To describe $d$, let $\varrho:\HZ_{(p)}\to \HF_p$ be reduction modulo $p$ and
$\partial:\HF_p\to\Sigma \HZ_{(p)}$ the connecting morphism for
multiplication by $p$. Put $\beta_{\HF_p}=(\Sigma\varrho)\partial$.\footnote{If
$a$ lifts a mod-$p$ cocycle $\bar a$ and $da=pb$, our sign convention is
$\partial[\bar a]=[b]$ and $\beta_{\HF_p}[\bar a]=[b\bmod p]$.} Set
\[
 j=\id_{\HF_p}\wedge\varrho:\HF_p\wedge \HZ_{(p)}\longrightarrow \HF_p\wedge \HF_p.
\]
The coefficient action on $\HF_p\wedge \HZ_{(p)}$ is through
$\id_{\HF_p}\wedge u:\HF_p\to \HF_p\wedge \HZ_{(p)}$. Its Bockstein is
$\beta_{\HF_p}\wedge\id_{\HZ_{(p)}}$; it preserves the splitting
$\HF_p\wedge \HZ_{(p)}\simeq \HF_p\oplus C/p$ and induces the differential $d$
on the second summand. We also write $d$ for the full Bockstein
on $\pi_{**}(\HF_p\wedge \HZ_{(p)})$, and denote its action on
$H=\pi_{**}\HF_p$ by $d_H$.

To fix notation for the following lemmas, we recall the Bockstein calculation of
Dundas and \O{}stv{\ae}r~\cite[Remark~1.1, Lemma~1.2 and formulas~(2)--(3)]{DO}
together with the corresponding decomposition of $\HF_p$-modules due to
Hoyois~\cite[Theorem~5.17]{Hoyois}.

\begin{proposition}\label{prop:bockstein-basis}
For every prime $p\ne\operatorname{char}k$, there are classes
$x_i,y_i\in\pi_{**}(\HF_p\wedge \HZ_{(p)})$, $i\ge1$, with
\begin{equation}\label{cc:eq:degrees}
 |x_i|=(2p^i-2,p^i-1),\qquad
 |y_i|=(2p^i-1,p^i-1),
\end{equation}
such that the ordered monomials $m=\prod_{i\ge1}x_i^{a_i}y_i^{\epsilon_i}$,
with $a_i\ge0$, $\epsilon_i\in\{0,1\}$ and finite support, form a basis
over the left coefficients $H$.
In particular, these classes define an isomorphism of left $\HF_p$-modules
in the homotopy category
\[
 \bigoplus_m\Sigma^{|m|}\HF_p\xrightarrow{\sim}\HF_p\wedge \HZ_{(p)}.
\]
Let $V$ be the $\Fp$-span of these ordered monomials. Write
$x^a=\prod_i x_i^{a_i}$ and let $i_1<\cdots<i_r$ be the indices with
$\epsilon_i=1$. Define $d_V$ on this basis by
\[
 d_V(x^a y_{i_1}\cdots y_{i_r})
 =\sum_{s=1}^r(-1)^{s-1}x^a x_{i_s}\,
 y_{i_1}\cdots\widehat{y_{i_s}}\cdots y_{i_r},
\]
and extend $\Fp$-linearly. Under the identification
$H\otimes_{\Fp}V\cong\pi_{**}(\HF_p\wedge\HZ_{(p)})$ given by
$c\otimes m\mapsto cm$, the Bockstein $d$ is the tensor product differential
\begin{equation}\label{cc:eq:total}
 d(c\otimes m)=d_H(c)\otimes m+(-1)^{|c|}c\otimes d_V(m),
\end{equation}
where $|c|$ is the first homological degree of the homogeneous coefficient $c$.
At $p=2$, let $\tau\in H_{0,-1}$ and $\rho\in H_{-1,-1}$ be the
classes represented by $-1$. Then
\begin{equation}\label{eq:bockstein-basis-two}
 y_i^2=\tau x_{i+1}+\rho y_{i+1},\qquad
 d_H(\tau)=\rho,\quad d_H(\rho)=0.
\end{equation}
\end{proposition}

\begin{proof}
Let $u_{\HF_p}=\varrho\circ u\colon\one\to \HF_p$ be the unit of $\HF_p$.
We use the left and right inclusions $\HF_p\to \HF_p\wedge \HF_p$ of
\cite[\S0.0.1]{DO}, given by $\eta_L=\id_{\HF_p}\wedge u_{\HF_p}$ and
$\eta_R=u_{\HF_p}\wedge\id_{\HF_p}$.
Let $\tau_i,\xi_i\in\pi_{**}(\HF_p\wedge \HF_p)$ be the generators defined in
\cite[\S12, p.~45, before Lemma~12.3]{VoevodskyOperations}, with the
conventions of Dundas and \O{}stv{\ae}r~\cite[p.~4, before Remark~1.1]{DO}.
Their Lemma~1.2 first describes the image of
$\HZ_{(p)}\wedge \HF_p\to \HF_p\wedge \HF_p$ as the right-coefficient algebra on
$\chi\tau_i,\chi\xi_i$ for $i\ge1$, where $\chi$ interchanges
the smash factors. Applying $\chi$ identifies the image of $j$
as the left-coefficient algebra on the nonconjugate
$\tau_i,\xi_i$, with the displayed monomial basis.
Since $j_*$ is injective, we define $x_i$ and $y_i$ by
$j_*(x_i)=\xi_i$ and $j_*(y_i)=\tau_i$.
The spectrum $\HF_p\wedge \HZ_{(p)}$ is a cellular $\HF_p$-module by
\cite[Theorem~5.17]{Hoyois}. Consequently the map from the free
module on this homotopy basis is an equivalence.

The Bockstein in \cite[Remark~1.1 and the paragraph following
Lemma~1.2]{DO} is explicitly that of the first smash factor;
it sends $\tau_i$ to $\xi_i$ and $\xi_i$ to zero. Note that
$d_V(y_i)=x_i$ and $d_V(x^a)=0$ by definition. Naturality
along $j$, its restriction $d_H$ on left coefficients, and the
graded Leibniz rule give \eqref{cc:eq:total}.
At $p=2$ the formula of \cite[Lemma~1.2]{DO}, after the same
interchange of factors, is precisely
\[
 \tau_i^2=\eta_L(\tau)\xi_{i+1}
              +\eta_L(\rho)\tau_{i+1}.
\]
This gives \eqref{eq:bockstein-basis-two}. Its compatibility
with the differential is explicit: applying $d$ to the
right-hand side gives
$\rho x_{i+1}+\rho x_{i+1}=0$.
\end{proof}

Let $\bar V\subset V$ be the span of the nonconstant monomials. These have
positive second homological degree, whereas
$\pi_{a,b}\HF_p=H^{-a}(k,\Fp(-b))=0$ for $b>0$. Thus their
augmentation is zero, the constant monomial gives the summand
$\HF_p$, and
\[
 \pi_{**}(C/p)=H\otimes_{\Fp}\bar V.
\]
The complex $(V,d_V)$ is the tensor product of the Koszul
complexes $\Fp[x_i]\otimes\Lambda(y_i)$. This is a statement
about its additive basis and differential; at $p=2$ its
multiplication is instead given by
\eqref{eq:bockstein-basis-two}.

\begin{lemma}\label{cc:lem:pairs}
There is a decomposition
\begin{equation}\label{cc:eq:pairs}
 \bar V=\bigoplus_\lambda(\Fp v_\lambda\oplus\Fp u_\lambda),\qquad
 d_Vv_\lambda=u_\lambda,
\end{equation}
with homogeneous $v_\lambda,u_\lambda$ of bidegrees
\begin{equation}\label{cc:eq:pair-degrees}
 |v_\lambda|=(2w_\lambda+\ell_\lambda,\,w_\lambda),\qquad
 |u_\lambda|=(2w_\lambda+\ell_\lambda-1,\,w_\lambda),\qquad
 w_\lambda\geq p-1,\quad\ell_\lambda\geq1.
\end{equation}
For each $w$ there are only finitely many $\lambda$ with $w_\lambda=w$.
\end{lemma}

\begin{proof}
A monomial $m$ has bidegree $(2w(m)+\ell(m),w(m))$, where
$w(m)=\sum_i(a_i+\epsilon_i)(p^i-1)$ is its weight and
$\ell(m)=\sum_i\epsilon_i$ its $y$-degree. So the bidegree of a
homogeneous element determines its weight and $y$-degree, and $d_V$
preserves the weight and lowers the $y$-degree by one. As a complex,
$V$ is the tensor product over $i\geq1$ of the Koszul complexes
$\Fp[x_i]\otimes\Lambda(y_i)$, each with homology $\Fp\cdot1$; by
the K\"unneth formula $H_*(V,d_V)=\Fp\cdot1$, so $\bar V$ is acyclic.
In each bidegree of $\bar V$, choose a complement to the cycles.
The differential $d_V$ maps this complement isomorphically onto the cycles
in the next bidegree. Bases $v_\lambda$
of these complements and $u_\lambda=d_Vv_\lambda$ give
\eqref{cc:eq:pairs}. Since $v_\lambda$ is not a cycle, $\ell_\lambda\geq1$;
nonconstant monomials have weight $\geq p-1$; and there are finitely
many monomials of each weight.
\end{proof}

\begin{lemma}\label{cc:lem:pair-cycles}
For each $\lambda$, $P_\lambda=H\otimes(\Fp v_\lambda\oplus\Fp u_\lambda)$
is a subcomplex of $(H\otimes\bar V,d)$, and $c\mapsto d(c\,v_\lambda)$
is a bijection from $H$ onto the cycles of $P_\lambda$ \textup{(}up to a
shift of bidegree\textup{)}. In particular $P_\lambda$, and hence
$H\otimes\bar V=\bigoplus_\lambda P_\lambda$, is acyclic.
\end{lemma}

\begin{proof}
By \eqref{cc:eq:total}, $d(a v_\lambda+b u_\lambda)
=d_H(a)\,v_\lambda+\bigl((-1)^{|a|}a+d_H(b)\bigr)u_\lambda$, which lies
in $P_\lambda$. This is zero if and only if
$a=-(-1)^{|a|}d_H(b)$, since then $d_H(a)=\pm d_H^2(b)=0$
automatically. Thus a cycle is determined by its arbitrary
$u_\lambda$-coefficient $b$. The $u_\lambda$-coefficient of
$d(c v_\lambda)$ is $(-1)^{|c|}c$, which runs through $H$ exactly once.
\end{proof}

\begin{lemma}\label{cc:lem:coop}
$p\,\pi_{**}C=0$, and $r\colon\pi_{**}C\to\pi_{**}(C/p)$ is injective
with image $\ker d=\im d$.
\end{lemma}

\begin{proof}
Write $\Delta=\pi_{**}C$, a $\Zp$-module, and
$\Delta[p]$ for its $p$-torsion. The acyclicity established in
Lemma~\ref{cc:lem:pair-cycles} gives $\ker d=\im d$. By
\eqref{cc:eq:bockstein-triangle}, $r$ induces an injection
$\Delta/p\Delta\hookrightarrow\pi_{**}(C/p)$ and $\delta$ maps
$\pi_{**}(C/p)$ onto $\Delta[p]$ (up to shift). Hence
\[
 \im d=r(\Delta[p]),\qquad
 \ker d=\delta^{-1}\bigl(\Delta[p]\cap p\Delta\bigr)\supseteq
 r(\Delta).
\]
As $\ker d=\im d\subseteq r(\Delta)$, we get
$\Delta=\Delta[p]+p\Delta$ and $\Delta[p]\cap p\Delta=0$.
Therefore $p\Delta=p^2\Delta$ has no $p$-torsion, so it is a
$\Q$-vector space. On the other hand $C_\Q=0$, because
$\HQ\wedge\HQ\simeq\HQ$, so $\Delta\otimes\Q=\pi_{**}(C_\Q)=0$ and
$\Delta$ is torsion. Hence $p\Delta=0$. Then $r$ is injective and
$r(\Delta)=\ker\delta=\ker d=\im d$.
\end{proof}

At the level of homotopy groups, the preceding lemma recovers
\cite[Theorem~1.4 and Corollary~1.5]{DO}. The following construction
realizes the calculation by actual maps of modules.

\begin{proposition}\label{cc:prop:modules}
There is an equivalence of left $\HZ_{(p)}$-modules
\begin{equation}\label{cc:eq:module-split}
 \HZ_{(p)}\wedge \HZ_{(p)}\simeq \HZ_{(p)}\oplus\bigoplus_\lambda
 \Sigma^{2w_\lambda+\ell_\lambda-1,\,w_\lambda}\HF_p,
\end{equation}
with $w_\lambda,\ell_\lambda$ as in \eqref{cc:eq:pair-degrees}.
\end{proposition}

\begin{proof}
Since $\Sigma^{|v_\lambda|}\HF_p$ is a free $\HF_p$-module on one generator,
$v_\lambda\in\pi_{|v_\lambda|}(C/p)$ is represented by an $\HF_p$-module map
$\widetilde v_\lambda\colon\Sigma^{|v_\lambda|}\HF_p\to C/p$. Define the
$\HZ_{(p)}$-module map
\begin{equation}\label{cc:eq:g}
 g_\lambda\colon\Sigma^{|u_\lambda|}\HF_p
 \xrightarrow{\ \Sigma^{-1}\widetilde v_\lambda\ }
 \Sigma^{-1}(C/p)\xrightarrow{\ \Sigma^{-1}\delta\ }C,
\end{equation}
using $|u_\lambda|=|v_\lambda|-(1,0)$. The first map is
$\HF_p$-linear, hence $\HZ_{(p)}$-linear, and the second is a map of $\HZ_{(p)}$-modules.

Identify $\pi_{**}(\Sigma^{|u_\lambda|}\HF_p)$ with $H$ (shifted). Since
$\widetilde v_\lambda$ is $H$-linear on homotopy,
$r\,g_\lambda(c)=\pm\,d(c\,v_\lambda)$, with a sign depending only on
the bidegree of $c$. By Lemma~\ref{cc:lem:pair-cycles}, the map
\[
 \textstyle\bigoplus_\lambda H\longrightarrow\pi_{**}(C/p),\qquad
 (c_\lambda)\longmapsto\sum_\lambda r\,g_\lambda(c_\lambda),
\]
is injective with image $\bigoplus_\lambda\ker(d|_{P_\lambda})=\ker d$,
which is $r(\pi_{**}C)$ by Lemma~\ref{cc:lem:coop}. Since $r$ is injective,
\[
 g=\textstyle\sum_\lambda g_\lambda\colon
 \bigoplus_\lambda\Sigma^{|u_\lambda|}\HF_p\longrightarrow C
\]
is an isomorphism on $\pi_{**}$ (homotopy commutes with sums because
$\Sigma^{a,b}\one$ is compact). Both sides are cellular $\HZ_{(p)}$-modules,
i.e.\ lie in the localizing subcategory generated by the $\Sigma^{a,b}\HZ_{(p)}$:
$\HF_p$ is the cofiber of $p$ on $\HZ_{(p)}$; $\HZ_{(p)}\wedge \HZ_{(p)}$ is the free module on the
cellular spectrum $\HZ_{(p)}$ \cite[Proposition~8.1]{Hoyois}; and $C$ is a
summand of $\HZ_{(p)}\wedge \HZ_{(p)}$. The cofiber of $g$ is a cellular $\HZ_{(p)}$-module with zero $\pi_{**}$.
Since $[\Sigma^{a,b}\HZ_{(p)},-]_{\HZ_{(p)}}=\pi_{a,b}(-)$ and these free modules
generate the cellular subcategory of $\Mod_{\HZ_{(p)}}$, the cofiber is zero.
Thus $g$ is an equivalence of $\HZ_{(p)}$-modules. Adding the summand $\HZ_{(p)}$
of \eqref{cc:eq:C} proves \eqref{cc:eq:module-split}.
\end{proof}

The following corollary is the main result of this subsection. It shows that $p$
annihilates the obstruction groups that occur in the proof of
Proposition~\ref{cc:prop:prime}.

\begin{corollary}\label{cc:cor:odd}
For all $a,b\in\Z$ there is a noncanonical isomorphism
\begin{equation}\label{cc:eq:operation-general}
 [\HZ_{(p)},\Sigma^{a,b}\HZ_{(p)}]\cong H^a(k,\Zp(b))\oplus
 \prod_\lambda H^{a-2w_\lambda-\ell_\lambda}(k,\Fp(b-w_\lambda)).
\end{equation}
Only the finitely many indices with $w_\lambda\le b$ contribute
in a fixed bidegree. In particular, for every $h\geq0$,
\begin{equation}\label{cc:eq:odd}
 p\,[\HZ_{(p)},\Sigma^{2h+1,h}\HZ_{(p)}]=0,\qquad
 [\HZ_{(p)},\Sigma^{2h+1,h}\HZ_{(p)}]=0\quad\text{if }h<p-1.
\end{equation}
\end{corollary}

\begin{proof}
By the free-forgetful adjunction, $[\HZ_{(p)},\Sigma^{a,b}\HZ_{(p)}]$ is the group of
$\HZ_{(p)}$-module maps $\HZ_{(p)}\wedge \HZ_{(p)}\to\Sigma^{a,b}\HZ_{(p)}$, which
Proposition~\ref{cc:prop:modules} computes as a product. The summand $\HZ_{(p)}$
contributes $\pi_{-a,-b}\HZ_{(p)}=H^a(k,\Zp(b))$. Applying the internal $\HZ_{(p)}$-module Hom to
$\HZ_{(p)}\xrightarrow{p}\HZ_{(p)}\to \HF_p$ gives
$\iHom_{\HZ_{(p)}}(\HF_p,\HZ_{(p)})\simeq\fib(p:\HZ_{(p)}\to \HZ_{(p)})\simeq\Sigma^{-1}\HF_p$,
hence
\[
 [\Sigma^{|u_\lambda|}\HF_p,\Sigma^{a,b}\HZ_{(p)}]_{\HZ_{(p)}}
 =\pi_{0,0}\,\Sigma^{(a,b)-|u_\lambda|-(1,0)}\HF_p
 =H^{a-2w_\lambda-\ell_\lambda}(k,\Fp(b-w_\lambda)).
\]
For $(a,b)=(2h+1,h)$ the first summand is $H^{2h+1}(k,\Zp(h))=0$,
since the degree exceeds the weight, and the remaining factors are
$\Fp$-vector spaces of weight $h-w_\lambda$. If $h<p-1\leq w_\lambda$,
the weights are negative and the factors vanish.
\end{proof}

\begin{corollary}\label{cor:global-cooperations}
There is a noncanonical equivalence of $\HZ$-modules
\[
 \HZ\wedge\HZ\simeq\HZ\oplus
 \bigoplus_p\bigoplus_{\lambda\in I_p}
 \Sigma^{2w_\lambda+\ell_\lambda-1,w_\lambda}\HF_p,
\]
where $I_p$ is the set of pairs in Lemma~\ref{cc:lem:pairs} for
the prime $p$. In particular $w_\lambda\ge p-1$.
\end{corollary}

\begin{proof}
Let $C_{\Z}$ be the fiber of multiplication
$\HZ\wedge\HZ\to\HZ$. Since $(C_{\Z})_{\Q}=0$, tensoring the
coefficient sequence $\Z\to\Q\to\Q/\Z$ with $C_{\Z}$ gives
$\Sigma C_{\Z}\simeq C_{\Z}\otimes_{\HZ}\mathrm H(\Q/\Z)$.
The decomposition
$\Q/\Z=\bigoplus_p\Q/\Z_{(p)}$ and the coefficient sequences
$\Z_{(p)}\to\Q\to\Q/\Z_{(p)}$ therefore give
$C_{\Z}\simeq\bigoplus_p(C_{\Z})_{(p)}$ as $\HZ$-modules.
Apply Proposition~\ref{cc:prop:modules} to each summand and
add the unit summand $\HZ$.
\end{proof}

\begin{corollary}
Let\label{cor:global-odd} $\Lambda\subseteq\Q$ be a localization of $\Z$. Then
$[\mathrm H\Lambda,\Sigma^{1,0}\mathrm H\Lambda]=0$.
For every $h\geq 0$, the integer
\[
 \prod_{\substack{p\le h+1\\1/p\notin\Lambda}}p,
\]
where $p$ runs over primes, annihilates
$[\mathrm H\Lambda,\Sigma^{2h+1,h}\mathrm H\Lambda]$.
\end{corollary}
\begin{proof}
Extend the module decomposition of
Corollary~\ref{cor:global-cooperations} to $\Lambda$. The free summand
contributes $H^{2h+1}(k,\Lambda(h))=0$, and the summands for primes
invertible in $\Lambda$ vanish upon localization. For each remaining prime $p$,
the contributions are $\Fp$-vector spaces of coefficient weight $h-w_\lambda$.
They vanish unless $p-1\le w_\lambda\le h$, so the displayed product
annihilates the group.
\end{proof}

\begin{remark}\label{cc:rem:first-pair}
The only nonconstant basis monomials of weight $p-1$ are
$y_1$ and $x_1$. They form a single pair with $\ell=1$, so
\eqref{cc:eq:operation-general} gives
\[
 [\HZ_{(p)},\Sigma^{2p-1,p-1}\HZ_{(p)}]\cong\Fp.
\]
This group is generated by
\[
 \partial P^1\varrho:\HZ_{(p)}\longrightarrow\Sigma^{2p-1,p-1}\HZ_{(p)},
\]
where $P^1=\mathrm{Sq}^2$ at $p=2$.
Here is a direct verification. Put $w=p-1$ and
$f=P^1\varrho:\HZ_{(p)}\to\T^w \HF_p$. By Proposition~\ref{prop:bockstein-basis} and the
vanishing of field motivic cohomology, evaluation on $x_1$
identifies
\[
 [\HZ_{(p)},\T^w \HF_p]=[\HF_p\wedge \HZ_{(p)},\T^w \HF_p]_{\HF_p}\xrightarrow{\sim}\Fp.
\]
The operation $f$ is nonzero: on $\mathbf P^p$ it sends the
reduction of $c_1(\mathcal O(1))$ to its nonzero $p$th power
\cite[Lemma~9.8]{VoevodskyOperations}. Hence
$\langle f,x_1\rangle\ne0$.

Since $d(y_1)=x_1$, Lemma~\ref{cc:lem:coop} gives a class
$z\in\pi_{2w,w}C$ with $r(z)=x_1$ and $pz=0$.
If $\partial f=0$, the coefficient exact sequence would lift
$f$ to an operation $\widetilde f:\HZ_{(p)}\to\T^w \HZ_{(p)}$.
Its $\HZ_{(p)}$-linear adjoint, restricted to $C$, sends $z$ into
$\pi_{0,0}\HZ_{(p)}=\Zp$, which has no $p$-torsion. Its value on $z$
is therefore zero. Reduction would give
$\langle f,x_1\rangle=0$, a contradiction. Thus
$\partial P^1\varrho$ is nonzero and generates the group.
\end{remark}

\subsection{Adams summands and extension obstructions}\label{cc:sec:adams}

Given a prime $p$ and a weight $q$, we seek an exponent $h$ such that
$p^h$ times the identity on the top slice $B_{q,(p)}$ extends to a map
$W_{q,(p)}\to B_{q,(p)}$.

Recall that a triangular matrix $T$ with diagonal entries
$\lambda_0,\ldots,\lambda_q$ satisfies $\prod_{i=0}^q(T-\lambda_i I)=0$.
When dealing with the $K_0$ group, we used the identity
$\ch_i(\psi^a x)=a^i\ch_i(x)$ to distinguish the weights as $a$ varies
in the proof of Riemann--Roch away from the small
primes~\cite[Proposition~4.3 and~(5.1)]{NavarroGarmendiaAway}.
The following two lemmas translate these observations into our setting.

\begin{lemma}\label{cc:lem:polynomial}
Let $\Lambda\subseteq\Q$ be a subring and let $\mathcal C$ be a $\Lambda$-linear triangulated
category. Let $0=Y_{q+1}\xrightarrow{\phi_q}Y_q\to\cdots
\xrightarrow{\phi_0}Y_0=Y$ be a tower in $\mathcal C$ with cofiber
sequences $Y_{i+1}\xrightarrow{\phi_i}Y_i\xrightarrow{\pi_i}G_i\to
\Sigma Y_{i+1}$. Let $t_i\in[Y_i,Y_i]$ satisfy
\[
 t_i\phi_i=\phi_it_{i+1},\qquad \pi_it_i=\lambda_i\pi_i
 \qquad(\lambda_i\in \Lambda).
\]
Then $\prod_{i=0}^q(t_0-\lambda_i)=0$ in $[Y,Y]$.
\end{lemma}

\begin{proof}
Put $\Phi_j=\phi_0\cdots\phi_{j-1}\colon Y_j\to Y$ and
$P_j=\prod_{i<j}(t_0-\lambda_i)$. We show by induction that
$P_j=\Phi_j\sigma_j$ for some $\sigma_j\colon Y\to Y_j$; for $j=0$ take
$\sigma_0=\id$. Since $t_0\Phi_j=\Phi_jt_j$, we have
$P_{j+1}=\Phi_j(t_j-\lambda_j)\sigma_j$. As $\pi_j(t_j-\lambda_j)=0$,
exactness gives $t_j-\lambda_j=\phi_j\tau$ for some
$\tau\colon Y_j\to Y_{j+1}$, so $\sigma_{j+1}=\tau\sigma_j$ works.
For $j=q+1$ we get $P_{q+1}=0$, since $Y_{q+1}=0$.
\end{proof}

\begin{lemma}\label{cc:lem:adams}
Let $\Lambda\subseteq\Q$ and $u\geq1$ with $1/u\in \Lambda$.
Then the $u$-th Adams operation $\Psi^u$ acts on the filtered object $f_\bullet W_{q,\Lambda}$ as in
Lemma~\ref{cc:lem:polynomial}, with $G_i=B_{i,\Lambda}$ and $\lambda_i=u^i$.
In particular $\prod_{i=0}^q(\Psi^u-u^i)=0$ in $[W_{q,\Lambda},W_{q,\Lambda}]$.
\end{lemma}

\begin{proof}
By functoriality, $t_i=f_i(\Psi^u)$ commutes with the inclusions and
$\pi_it_i=s_i(\Psi^u)\pi_i$. Now
$s_i(\Psi^u)\in[B_{i,\Lambda},B_{i,\Lambda}]=\Lambda$ by Lemma~\ref{cc:lem:descent}.
Since $\Lambda\hookrightarrow\Q$, it suffices to compute the rational scalar. The weight-$i$ character satisfies
$\ch_i\circ\Psi^u=u^i\ch_i$ by the action of Adams operations on the
rational weight summands \cite[Proposition~5.3.14 and Remark~6.2.3.10]{Riou}.
Taking $s_i$ gives the scalar $u^i$ since $s_i(\ch_i)$ is an equivalence
by Lemma~\ref{lem:slice-normalization}.
\end{proof}

Fix a prime $p$ and let $q\geq0$ be the top weight of the truncation
$W_{q,(p)}$. Choose an integer $u>1$ whose residue generates $\Fp^\times$,
and write $[j]$ for the class of $j$ modulo $p-1$.
For each class $a\in\Z/(p-1)$ represented by a weight $0\leq i\leq q$, put
\[
 P_a(t)=\prod_{\substack{0\leq i\leq q\\ i\equiv a\ (p-1)}}(t-u^i)
 \in\Zp[t].
\]
If $i\not\equiv j\pmod{p-1}$, then $u^i\not\equiv u^j\pmod p$, so
$u^i-u^j\in\Zp^\times$. For $a\neq b$, the resultant $\Res(P_a,P_b)$
is a product of these units. Thus $P_a$ and $P_b$ are comaximal, and the Chinese remainder theorem
gives orthogonal idempotents $e_a\in\Zp[t]/(\prod_aP_a)$ with
$\sum_ae_a=1$ and $e_a\equiv\delta_{ab}\pmod{P_b}$. By
Lemma~\ref{cc:lem:adams}, $t\mapsto\Psi^u$ defines a ring map
$\Zp[t]/(\prod_aP_a)\to[W_{q,(p)},W_{q,(p)}]$, so the $e_a(\Psi^u)$ are
orthogonal idempotents with sum $\id$. As $u^j$ is a root of $P_{[j]}$,
\begin{equation}\label{cc:eq:idempotent-slices}
 s_j\bigl(e_a(\Psi^u)\bigr)=e_a(u^j)=
 \begin{cases}1&j\equiv a\pmod{p-1},\\0&\text{otherwise.}\end{cases}
\end{equation}
For $p=2$ all weights lie in a single residue class, so this decomposition
is unnecessary and $e_{[q]}(\Psi^u)=\id$.

\begin{proposition}\label{cc:prop:prime}
Let $\iota_q\colon B_q\to W_q$ denote the inclusion of the top slice and
put $h=h_p(q)=\lfloor q/(p-1)\rfloor$. There is a map
\begin{equation}\label{cc:eq:nu}
 \nu_{q,p}\colon W_{q,(p)}\longrightarrow B_{q,(p)},\qquad
 \nu_{q,p}\,\iota_{q,(p)}=p^{h}\id.
\end{equation}
In particular, for $q\geq1$ the connecting morphism
$\kappa_q\colon W_{q-1}\to\Sigma B_q$ in \eqref{cc:eq:attaching}
satisfies $p^{h}\kappa_{q,(p)}=0$.
\end{proposition}

\begin{proof}
Let $e=e_{[q]}(\Psi^u)$ and split it in the idempotent-complete category
$\SH(k)_{(p)}$: $W_{q,(p)}\xrightarrow{r}E\xrightarrow{i}W_{q,(p)}$ with
$ri=\id_E$ and $ir=e$. The exact functor $s_j$ sends this retraction to
the splitting of $s_j(e)$, so by \eqref{cc:eq:idempotent-slices}, $s_jE=0$
unless $j$ lies in $J=\{q-t(p-1):0\leq t\leq h\}$, and $s_j(r)$ is an
equivalence $B_{j,(p)}\to s_jE$ for $j\in J$. Also $f_{q+1}E=0$, as a
retract of $f_{q+1}W_{q,(p)}=0$.

Put $E^t=f_{q-t(p-1)}E$ for $0\leq t\leq h$, and let
$j_t\colon E^t\to E$ be the canonical map. Since $f_{j+1}E\to f_jE$ is
an equivalence for $j\notin J$, $j_h$ is an equivalence,
$E^0=f_qE\simeq s_qE$, and there are cofiber sequences
\begin{equation}\label{cc:eq:block-tower}
 E^{t-1}\longrightarrow E^t\longrightarrow B_{q-t(p-1),(p)}
 \xrightarrow{\ \partial_t\ }\Sigma E^{t-1}\qquad(1\leq t\leq h).
\end{equation}
The map $f_q(r)\colon B_{q,(p)}=f_qW_{q,(p)}\to E^0$ equals $s_q(r)$, so it
is an equivalence; and by naturality of $f_q\to\id$,
\begin{equation}\label{cc:eq:top-restriction}
 r\,\iota_{q,(p)}=j_0\circ f_q(r).
\end{equation}

We construct $v_t\colon E^t\to B_{q,(p)}$ with $v_0=f_q(r)^{-1}$ and
$v_t|_{E^{t-1}}=p\,v_{t-1}$. Given $v_{t-1}$, the obstruction to
extending a map $x\colon E^{t-1}\to B_{q,(p)}$ across
\eqref{cc:eq:block-tower} is
\[
 x\circ\Sigma^{-1}\partial_t\in[\Sigma^{-1}B_{q-t(p-1),(p)},B_{q,(p)}]
 \cong[\HZ_{(p)},\Sigma^{2n+1,n}\HZ_{(p)}],\qquad n=t(p-1)>0,
\]
since $B_{j,(p)}\simeq\T^j \HZ_{(p)}$. Corollary~\ref{cc:cor:odd} shows that this group is killed by
$p$, so $x=p\,v_{t-1}$ extends to some $v_t$. Now set
$\nu_{q,p}=v_h\circ j_h^{-1}\circ r$. By \eqref{cc:eq:top-restriction},
\[
 \nu_{q,p}\,\iota_{q,(p)}
 =v_h\,j_h^{-1}j_0\,f_q(r)=(v_h|_{E^0})\,f_q(r)=p^h v_0\,f_q(r)=p^h.
\]
The last assertion follows from the $p$-local form of
\eqref{cc:eq:criterion}.
\end{proof}

Without the Adams decomposition, extending across each individual
slice would give $p^{\max(0,q-p+2)}$: only the gaps below $p-1$
have zero obstruction group. The Adams summand reduces this to
$p^{\lfloor q/(p-1)\rfloor}$, using the exponent-$p$ statement
of Corollary~\ref{cc:cor:odd} at every remaining step.

\subsection{The stable character by arithmetic gluing}\label{cc:sec:glue}

Put $\Lambda_q=\Z[1/T_q]$. A finite difference of Adams operations gives the
top projection over $\Lambda_q$.

\begin{proposition}\label{cc:prop:D}
For $q\geq0$ there is a map
\[
 \theta_q^{\mathrm{loc}}\colon W_{q,\Lambda_q}\longrightarrow B_{q,\Lambda_q},\qquad
 \theta_q^{\mathrm{loc}}\,\iota_{q,\Lambda_q}=\id.
\]
For $q>0$ it is obtained from the endomorphism
\begin{equation}\label{cc:eq:finite-difference}
 \varepsilon=\frac1{q!}\sum_{v=0}^{q}(-1)^{q-v}\binom qv\Psi^{v+1}
 \quad\text{of }W_{q,\Lambda_q}.
\end{equation}
\end{proposition}

\begin{proof}
For $q=0$ take $\theta_0^{\mathrm{loc}}=\id$, as $W_0=B_0$. Let $q>0$. The integers
$1,\ldots,q+1$ and $q!$ are units in $\Lambda_q$ by
Lemma~\ref{lem:denominators}, so \eqref{cc:eq:finite-difference} is
defined. By Lemma~\ref{cc:lem:adams}, $\varepsilon$ acts on the filtered object
$f_\bullet W_{q,\Lambda_q}$ as in Lemma~\ref{cc:lem:polynomial}, with
\[
 \lambda_i=\frac1{q!}\sum_{v=0}^q(-1)^{q-v}\binom qv(v+1)^i
 =\frac{(\Delta_+^qx^i)(1)}{q!},
\]
where $(\Delta_+ f)(x)=f(x+1)-f(x)$ is the forward difference. Since $\Delta_+^q$ kills
polynomials of degree $<q$ and $\Delta_+^qx^q=q!$, we get
\[ \lambda_i=\delta_{iq}\qquad(0\le i\le q). \] Lemma~\ref{cc:lem:polynomial} gives
$\varepsilon^q(\varepsilon-1)=0$, so $\varepsilon^{q+1}=\varepsilon^q$ and $P=\varepsilon^q$ is idempotent.

Split $P$ as $W_{q,\Lambda_q}\xrightarrow{r}E\xrightarrow{i}W_{q,\Lambda_q}$. Then
$s_jE$ is the image of $s_j(P)=\lambda_j^q$, so $s_jE=0$ for $j<q$,
$s_q(r)$ is an equivalence, and $f_{q+1}E=0$. Hence $f_qE\to E$ is an
equivalence, and so is $r\,\iota_{q,\Lambda_q}=(f_qE\to E)\circ f_q(r)$. Put
$\theta_q^{\mathrm{loc}}=(r\,\iota_{q,\Lambda_q})^{-1}r$.
\end{proof}

\begin{remark}\label{rem:alternate-local-projection}
Alternatively, Corollary~\ref{cor:global-cooperations} shows that
$[\HZ_{\Lambda_q},\Sigma^{2h+1,h}\HZ_{\Lambda_q}]=0$ for
$0\le h\le q$: the surviving primes satisfy $p>q+1$ and their
weights exceed $q$. The top identity therefore extends across
all lower slices with zero obstruction.
\end{remark}

We use a finite arithmetic square for motivic cohomology coefficients.
Its classical stable analogue is Bousfield's arithmetic square for
localization at complementary sets of primes
\cite[Proposition~2.10]{Bousfield}.

\begin{lemma}\label{cc:lem:arithmetic}
Let $S$ be a nonempty finite set of primes and put $\Lambda=\Z[S^{-1}]$.
The square of localization maps
\begin{equation}\label{cc:eq:arithmetic}
 \begin{tikzcd}[column sep=large,row sep=large]
 \HZ\arrow[r]\arrow[d]&\mathrm H\Lambda\arrow[d]\\
 \displaystyle\bigoplus_{p\in S}\mathrm H\Zp\arrow[r]&
 \displaystyle\bigoplus_{p\in S}\HQ
 \end{tikzcd}
\end{equation}
is a pullback in $\SH(k)$.
\end{lemma}

\begin{proof}
The sequence
\begin{equation}\label{cc:eq:arithmetic-exact}
 0\longrightarrow\Z\longrightarrow\Lambda\oplus\bigoplus_{p\in S}\Zp
 \xrightarrow{\ (a,(b_p))\mapsto(a-b_p)_p\ }\bigoplus_{p\in S}\Q
 \longrightarrow0
\end{equation}
is exact. The kernel of the difference map is the diagonal copy of
$\Lambda\cap\bigcap_{p\in S}\Zp=\Z$. Given
$(x_p)$, choose $y_p\in\Z[1/p]$ with $x_p-y_p\in\Zp$ (the principal part
of $x_p$ at $p$). Then $a=\sum_py_p\in\Lambda$ and $a-x_p\in\Zp$ for each $p$.
The exact coefficient sequence gives a distinguished triangle of
motivic Eilenberg--Mac Lane spectra. In a stable category the
associated square is both a pushout and a pullback, proving
\eqref{cc:eq:arithmetic}.
\end{proof}

\goodbreak
We now prove Theorem~\ref{cc:thm:main}, stated at the beginning of this section.

\begin{integralrealizationagain}[Integral realization]
\IntegralRealizationStatement{equation*}{\tag{\ref*{cc:eq:Phi-intro}}}
\end{integralrealizationagain}

\begin{proof}
(a) For $q=0$ take $\theta_0=\id$. Let $q>0$, let $S$ be the set of
primes $p\leq q+1$, so that $\Lambda_q=\Z[S^{-1}]$. For each $p\in S$,
write $h_p=h_p(q)=\lfloor q/(p-1)\rfloor$ for the exponent of $p$ in $T_q$.
Precomposing with localization maps, consider
\[
 a=T_q\,\theta_q^{\mathrm{loc}}\in[W_q,B_{q,\Lambda_q}],\qquad
 b_p=\frac{T_q}{p^{h_p}}\,\nu_{q,p}\in[W_q,B_{q,(p)}]\quad(p\in S),
\]
from Propositions~\ref{cc:prop:D} and~\ref{cc:prop:prime}. All have top value
$T_q$, so by Lemma~\ref{cc:lem:rational} they have the same image $T_q\ch_q$
in $[W_q,B_{q,\Q}]=[W_{q,\Q},B_{q,\Q}]$. Applying $[W_q,-]$ to the
$\T^q$-twist of \eqref{cc:eq:arithmetic} gives an exact Mayer--Vietoris
sequence
\[
 [W_q,B_q]\longrightarrow[W_q,B_{q,\Lambda_q}]\oplus\bigoplus_{p\in S}
 [W_q,B_{q,(p)}]\longrightarrow\bigoplus_{p\in S}[W_q,B_{q,\Q}],
\]
whose second map is $(a,(b_p))\mapsto(a_\Q-b_{p,\Q})_p$. So
$(a,(b_p))$ lifts to some $\theta_q\in[W_q,B_q]$. Its top value lies in
$[B_q,B_q]=\Z$ and maps to $T_q\in \Lambda_q$, hence equals $T_q$. The
equivalence with $T_q\kappa_q=0$ is \eqref{cc:eq:criterion}. The gluing is performed on a finite pullback of coefficient
spectra; applying $\Map(W_q,-)$ preserves this pullback.

(b) The first assertion is Lemma~\ref{cc:lem:rational}. If $\theta$ and
$\theta'$ have the same top value, exactness of $[-,B_q]$ on
\eqref{cc:eq:attaching} shows that $\theta-\theta'$ factors through
$W_q\to W_{q-1}$, and $[W_{q-1,\Q},B_{q,\Q}]=0$ by
Lemma~\ref{cc:lem:rational}.

(c) The coefficients $T_M/T_q$ are integers by
Lemma~\ref{lem:denominators}. Since $f_q(\theta_q)=\theta_q\iota_q$
and $s_q(\pi_{M,q})=\id$, we get
$s_q((\Phi_M)_q)=(T_M/T_q)T_q=T_M$, while $s_q$ of the other components
vanishes because $s_q(B_{q'})=0$ for $q'\neq q$. The rationalization is
given by (b). After inverting $T_M$, $\Phi_M$ is an isomorphism on all
slices; source and target have $f_{M+1}=0$, so its cofiber has a finite
slice filtration with zero slices and vanishes.
\end{proof}

For the fixed choices of $\theta_q$, the finite characters satisfy
\begin{equation}\label{eq:window-compatibility}
 \operatorname{pr}_{\le M}\Phi_{M'}
 =\frac{T_{M'}}{T_M}\Phi_M\pi_{M',M}\qquad(M'\ge M).
\end{equation}
Here $\operatorname{pr}_{\le M}\colon H_{M'}\to H_M$ is the projection onto
the summands $B_0,\ldots,B_M$. The identity follows from the componentwise
definition of the finite characters.

\subsection{Evaluations and denominators}\label{sec:thom}

In Subsection~\ref{cc:sec:glue} we constructed an integral map
$\Phi_M\colon W_M\to H_M$ for each $M\geq0$. We now identify the classes that its components induce in integral motivic
cohomology, both on $x\in K_s(X)$ and on its product with a Thom class.
Finally, we prove that $T_q$ is the least positive integer $N$ for which
$N\ch_q$ admits an integral stable realization.

Let $q_M\colon A_0\to W_M$ be the quotient. For $a\le b$, define the
translates of $\Phi_M$ by
\begin{equation}\label{eq:psi}
 \psi^{(M)}_{a,b}
 =\T^a(\Phi_Mq_M)\,\iota_{b,a}
 \colon A_b\longrightarrow\bigoplus_{t=0}^M B_{a+t}.
\end{equation}
Using the Bott equivalence $A_a\simeq\T^aA_0$, we can write
\eqref{eq:psi} as the composite
\[
\begin{tikzcd}[column sep=small]
 A_b \arrow[r,"\iota_{b,a}"] & A_a \arrow[r,"\sim"] &
 \T^aA_0 \arrow[r,"\T^aq_M" {midway}] &[1.6em] \T^aW_M
 \arrow[r,"\T^a\Phi_M" {midway}] &[1.6em] \displaystyle\bigoplus_{t=0}^M B_{a+t}.
\end{tikzcd}
\]
By Theorem~\ref{cc:thm:main}(c), the $t$-th component of
$(\Phi_M q_M)_\Q$ is the restriction of $T_M\ch_t$ to $A_{0,\Q}$.
To see what happens under the Bott equivalence, recall the normalization
$\ch_i\circ\beta^i_\Q=\T^i\ch_0$ used in the proof of
Lemma~\ref{lem:slice-normalization} \cite[Definition~6.2.3.9]{Riou}.
It gives
\[
 \ch_{a+t}\circ\beta^a_\Q=\T^a\ch_t
 \colon \T^a\KGL_\Q\longrightarrow \T^{a+t}\HQ.
\]
Indeed, precomposing both sides with the equivalence $\T^a\beta^t_\Q$
gives $\T^{a+t}\ch_0$, since $\beta^a\circ\T^a\beta^t=\beta^{a+t}$.
Thus transporting the $t$-th component from $\T^aA_0$ to $A_a$ gives
$T_M\ch_{a+t}$ on $A_{a,\Q}$. Precomposition with $\iota_{b,a}$ further
restricts this map along $A_b\to A_a\to\KGL$.

We work in the pointed unstable motivic homotopy category
$\mathcal H_\bullet(k)$. The $\mathbf P^1$-suspension spectrum functor
$\Sigma^\infty$ and the $\mathbf P^1$-infinite loop space functor
$\Omega^\infty$ form the adjoint pair
\[
 \Sigma^\infty\colon\mathcal H_\bullet(k)
 \rightleftarrows\SH(k):\Omega^\infty.
\]
Algebraic $K$-theory is represented by
$\cK=\Omega^\infty\kgl\simeq\Z\times\BGL$, pointed at the zero class,
where $\BGL\simeq\operatorname*{colim}_N\Gr(N,2N)$
\cite[\S4, Proposition~3.7 and Theorem~3.13]{MV}.
Let $\epsilon\colon\Sigma^\infty\cK\to A_0$ be the counit.
Effective-cover adjunction and representability identify maps from
suspension spectra to $A_0$ with reduced algebraic $K$-classes.
The ring structure described in Section~\ref{sec:preliminaries} gives the action
\[
 A_e\wedge A_0\longrightarrow A_e.
\]
For a vector bundle $V$ on $X$ with zero section $s$, its unstable Thom
space over $k$ is the pointed motivic space $\Thsp(V)=V/(V\setminus s(X))$.
For the universal rank-$e$ bundle $\gamma_e$ on
$\BGL_e\simeq\operatorname*{colim}_N\Gr(e,e+N)$, put
$U_e=\Thsp(\gamma_e)\wedge\cK$.

\begin{lemma}\label{lem:detection}
For $U=\cK$ or $U_e$, the coefficient map
$\widetilde H^{2w}(U,\Z(w))\to\widetilde H^{2w}(U,\Q(w))$
between the reduced motivic cohomology groups is injective.
The groups in degree $2w$ vanish for $w<0$ on $\cK$ and for $w<e$ on $U_e$.
\end{lemma}

\begin{proof}
On one rank component, put $Y=\Thsp(\gamma_e)\wedge\BGL_+$ and
$G_N=\Gr(e,e+N)\times\Gr(N,2N)$. The stabilization maps of Grassmannians
induce split injections on motives \cite[proof of Theorem~6.2.1.2]{Riou},
hence so do their products. Thus the restrictions in motivic cohomology
are surjective in every bidegree. By the Thom isomorphism at each finite
stage and the Milnor exact sequence \cite[pp.~221--222]{PaninPimenovRondigs},
\[
 \widetilde H^{2w}(Y,R(w))\cong
 \varprojlim_N\bigl(\CH^{w-e}(G_N)\otimes R\bigr)
 \qquad(R=\Z,\Q).
\]
Each Chow group is free on Schubert classes and vanishes in negative
codimension, proving both assertions for $Y$. Taking the product over ranks
gives the result for $\Thsp(\gamma_e)\wedge\cK_+$, where $+$ adds a disjoint
base point. The split inclusion $S^0\to\cK_+$ selecting
the zero class has cofiber $\cK$, so $U_e$ is a stable summand and inherits
both assertions. Taking $e=0$ gives $U_0=\cK$.
\end{proof}

The universal Newton class $s_q$ on $\mathcal K$ is characterized by
additivity and $s_q([L])=c_1(L)^q$ for line bundles; $s_0=\rk$.
It satisfies $s_q=q!\ch_q$ rationally. For $x\in K_n(X)$, pullback of
the universal class $c_q$ along a representative $S^n\wedge X_+\to\cK$,
followed by the suspension isomorphism, gives the higher Chern class
$c_{q,n}(x)$. This is Gillet's construction. For its motivic formulation
and normalization, see \cite[Remark~2.1.18]{DegliseOrientation} and
\cite{Pushin}. In particular no factorial is included in $c_{q,n}$.
For $q\geq1$, Newton's identity has the form
\[
 s_q=(-1)^{q-1}q\,c_q+P_q(c_1,\ldots,c_{q-1}),
\]
where each monomial of $P_q$ is a product of at least two Chern classes
of positive degree less than $q$. For $n>0$, cup products of reduced classes
on $S^n\wedge X_+$ vanish because its reduced diagonal is nullhomotopic.
Thus every monomial of $P_q$ pulls back to zero, giving
$\chi_{q,n}(x)=(-1)^{q-1}q\,c_{q,n}(x)$.

\begin{lemma}
Let\label{lem:evaluation} $a\le0$, $0\le t\le M$ and $q=a+t$. For any pointed motivic space
$U$ and any map $z\colon\Sigma^\infty U\to A_0$, one has
\[
 \psi^{(M)}_{a,0}(z)_t=
 \begin{cases}(T_M/q!)\,s_q(\bar z),&q\ge0,\\0,&q<0,
 \end{cases}
\]
where $\psi^{(M)}_{a,0}$ is the operation defined in \eqref{eq:psi}
and $\bar z\colon U\to\cK$ denotes the adjoint of $z$.
When $U=S^n\wedge X_+$ with $n\geq0$ and $z$ represents $x\in K_n(X)$,
the suspension isomorphism identifies $s_q(\bar z)$ with
$\chi_{q,n}(x)\in\CH^q(X,n)$ for $q\geq0$.
\end{lemma}

\begin{proof}
Let $\epsilon\colon\Sigma^\infty\cK\to A_0$ be the counit, adjoint to
$\id_{\cK}$. For $q\geq0$, the rationalization of \eqref{eq:psi} gives
\[
 \bigl(\psi^{(M)}_{a,0}(\epsilon)_t\bigr)_\Q
 =T_M\ch_q(\epsilon)
 =\frac{T_M}{q!}\,(s_q)_\Q.
\]
Since $q\leq M$, we have $q!\mid T_M$ by \eqref{eq:divisibility}, so
Lemma~\ref{lem:detection} gives the equality integrally. For $q<0$,
the same lemma gives $\widetilde H^{2q}(\cK,\Z(q))=0$.
Now $z=\epsilon\circ\Sigma^\infty\bar z$, hence
\[
 \psi^{(M)}_{a,0}(z)_t
 =\bar z^*\!\left(\psi^{(M)}_{a,0}(\epsilon)_t\right)
 =\begin{cases}
    (T_M/q!)\,s_q(\bar z),&q\geq0,\\
    0,&q<0.
  \end{cases}
\]
For $U=S^n\wedge X_+$ and $q\geq0$, the suspension isomorphism sends
$s_q(\bar z)$ to $\chi_{q,n}(x)$ by the Newton formula \eqref{eq:newton}.
\end{proof}

In particular, the weight-$q$ component of $\Phi_M$ on $K_n(X)$ is
$(T_M/q!)\chi_{q,n}$. These induced operations are independent of
the choices of the $\theta_q$.

For an oriented motivic ring spectrum $E$ and a virtual bundle $\xi$ of rank
$r$ on $X$, we denote by $\ThomClass{E}(\xi)\in E^{2r,r}(\Th(\xi))$ its Thom
class. On the trivial line it is the suspended unit. For a line bundle $L$
with zero section $s$, we have
$s^*\ThomClass{\KGL}(L)=c_1^{\KGL}(L)=\beta^{-1}(1-[L^\vee])$.
The rational Chern character satisfies
\begin{equation}\label{eq:rational-thom}
 \ch(\ThomClass{\KGL}(\xi))=\ThomClass{\HQ}(\xi)\Td(\xi)^{-1}.
\end{equation}
By multiplicativity of the Thom and Todd classes and the splitting principle,
this classical calculation reduces to a line bundle $L$, for which the
factor is $(1-e^{-x})/x=\Td(L)^{-1}$ with $x=c_1^{\HQ}(L)$
(compare \cite[Remark~5.1.6 and Theorem~6.3.1]{Riou}).

The lift of a Thom class to an effective cover is canonical. For a bundle
$V$ of rank $e$ on $X$, the spectrum $\Th(V)$ is $e$-effective.
Indeed, on a trivializing open $U$ its Thom spectrum over $k$ is
$\T^e\Sigma^\infty U_+$, and the assertion follows by descent.
Under Bott periodicity, $\ThomClass{\KGL}(V)$ gives a map $\Th(V)\to\KGL$.
Effective-cover adjunction lifts this map uniquely to $A_e\simeq\T^e\kgl$.
This lift is $\ThomClass{\kgl}(V)$ and is compatible with the action
$A_e\wedge A_0\to A_e$. After pullback to an affine bundle torsor,
write $\xi=V-\mathcal O^h$ with $e=r+h$. Desuspending by $\T^h$ gives
$\ThomClass{\kgl}(\xi)\in[\Th(\xi),A_r]$. Its independence of this presentation
and descent follow from the same adjunction and the corresponding
properties of $\ThomClass{\KGL}(\xi)$.

\begin{proposition}\label{prop:thom}
Let $\xi$ be a virtual bundle of rank $r$ on $X$, and let
$x\in K_n(X)$ with $n\geq0$. Put
\[
 z=\ThomClass{\kgl}(\xi)x
 \in[\Sigma^{n,0}\Th(\xi),A_r].
\]
For $a\le r$ and $0\le t\le M$, the Thom isomorphism identifies
\begin{equation}\label{eq:thom}
 \psi^{(M)}_{a,r}(z)_t
 =T_M[\ch(x)\Td(\xi)^{-1}]^{(n)}_{a+t-r}.
\end{equation}
Both sides are zero when $a+t-r<0$.
\end{proposition}

\begin{proof}
By Jouanolou's device \cite{Jouanolou}, replace $X$ by an affine
vector-bundle torsor. This is an $\mathbb A^1$-equivalence, also for
the pulled-back Thom spectra. On the affine scheme write
$\xi=V-\mathcal O^h$, with $V$ a bundle of rank $e=r+h$.
Suspending by $\T^h$ turns the assertion into the assertion for
$\psi^{(M)}_{a+h,e}$ applied to $\ThomClass{\kgl}(V)x$.

The universal class is
\[
 u=\ThomClass{\kgl}(\gamma_e)\boxtimes\epsilon:
 \Sigma^\infty U_e\longrightarrow A_e.
\]
Let $a_x:S^n\wedge X_+\to\mathcal K$ represent $x$, let
$t_V:\Thsp(V)\to\Thsp(\gamma_e)$ be induced by the classifying map of
$V$, and let $\delta_V:\Thsp(V)\to\Thsp(V)\wedge X_+$ be the Thom
diagonal. For every $n\ge0$, the required universal map is the composite
\[
 \begin{split}
 S^n\wedge\Thsp(V)
 &\xrightarrow{1\wedge\delta_V}
 S^n\wedge\Thsp(V)\wedge X_+\\
 &\simeq\Thsp(V)\wedge S^n\wedge X_+
 \xrightarrow{t_V\wedge a_x}\Thsp(\gamma_e)\wedge\mathcal K.
 \end{split}
\]
Its suspension pulls $u$ back to $\ThomClass{\kgl}(V)x$.

Put $j=a+t-r$ and $w=a+h+t=e+j$. If $j<0$, then $w<e$ and
Lemma~\ref{lem:detection} gives $\widetilde H^{2w}(U_e,\Z(w))=0$. Otherwise $j\le t\le M$, and the two integral
classes on the universal Thom space to be compared are
\[
 \psi^{(M)}_{a+h,e}(u)_t,
 \qquad
 \ThomClass{\HZ}(\gamma_e)\,
 T_M[\ch(\epsilon)\Td(\gamma_e)^{-1}]_j.
\]
The second is integral by \eqref{eq:integral-polynomial}. They agree
rationally by \eqref{eq:rational-thom} and the normalization of
\eqref{eq:psi}; Lemma~\ref{lem:detection} makes them equal integrally.
Pulling back proves the claim for $n=0$.
For $n>0$, the Thom diagonal makes this pullback linear in the
$\mathcal K$ variable: decomposable reduced classes vanish, while
classes from $BGL_e$ act through their pullbacks to $X$. Explicitly,
with $j=a+t-r$, the resulting class under the Thom isomorphism is
\[
 T_M\sum_{i=1}^j
 \frac{(-1)^{i-1}}{(i-1)!}
 [\Td(\xi)^{-1}]_{j-i}\,c_{i,n}(x),
\]
which is the right side of \eqref{eq:thom}.\qedhere
\end{proof}

\begin{theorem}[Exact stable denominator]\label{thm:exact-denominator}
For $q\ge0$, the image of the top restriction map
\[
 [W_q,B_q]\longrightarrow[B_q,B_q]=\Z
\]
is $T_q\Z$. Equivalently, the image of
$[\kgl,\T^q\HZ]$ under rationalization is
$\Z\cdot T_q\ch_q$. For $q\ge1$, the attaching map
$\kappa_q\colon W_{q-1}\longrightarrow\Sigma B_q$ has order exactly $T_q$.
\end{theorem}

\begin{proof}
Existence in Theorem~\ref{cc:thm:main} gives the inclusion
$T_q\Z$ in the image. Conversely, let $\theta\colon W_q\to B_q$
have top value $N$, and write
$\alpha\colon A_0\to W_q\xrightarrow{\theta}B_q$.
Its rationalization is $N\ch_q$ by Lemma~\ref{cc:lem:rational}.
For the universal rank-$b$ bundle $\gamma_b$ on $BGL_b$, evaluate
the translate $\T^b\alpha\colon A_b\to B_{b+q}$ on
$\ThomClass{\kgl}(\gamma_b)$. The integral Thom isomorphism
gives a class in $\CH^q(BGL_b)$ whose rational value, by
\eqref{eq:rational-thom}, is
\[
 N[\Td(\gamma_b)^{-1}]_q.
\]
Since $\CH^*(BGL_b)=\Z[c_1,\ldots,c_b]$, this is an integral
symmetric polynomial in the Chern roots for every $b$.

Fix a prime $p\le q+1$ and write
$q=h(p-1)+r$, where $0\le r<p-1$. In a sufficient number of
variables, the coefficient of the monomial symmetric function
of partition $((p-1)^h,r)$ (omitting the zero part when $r=0$)
in the degree-$q$ part of
\[
 \prod_{j=1}^b\frac{1-e^{-x_j}}{x_j}
 =\prod_{j=1}^b\sum_{a\ge0}\frac{(-1)^a x_j^a}{(a+1)!}
\]
is $(-1)^q/((p!)^h(r+1)!)$. Its $p$-adic valuation is $-h$.
Integrality therefore forces $v_p(N)\ge h$ for every such $p$,
and hence $T_q\mid N$. The descent isomorphism of
Lemma~\ref{cc:lem:descent} gives the second assertion, and the
extension criterion \eqref{cc:eq:criterion} gives the last.
\end{proof}

\begin{remark}\label{rem:exact-low-weights}
The attaching map $\kappa_1$ has order $2$, and
$(\kappa_1)_{(2)}=\partial\mathrm{Sq}^2\varrho$ in
$[\mathrm H\Z_{(2)},\Sigma^{3,1}\mathrm H\Z_{(2)}]$.
Its mod-$2$ reduction is the restriction of the first Milnor operation,
in agreement with the first slice differential of $\KGL/2$
computed in \cite[Lemma~5.1]{RondigsOstvaer}.
Moreover
$\kappa_2$ has order $12=2^2\cdot3$. At $p=2$ the construction
of the second top projection has two extension steps, and at
$p=3$ it has one. These bounds are forced simultaneously by
the universal Thom calculation.
\end{remark}

\section{Integral Riemann--Roch formula}\label{sec:transfer}

Let $k$ be a field of characteristic zero, let $X$ and $Y$ be smooth
quasi-projective $k$-schemes of pure dimension, and let $f\colon X\to Y$
be a projective $k$-morphism. Write $p_X$ and $p_Y$ for their structure
maps to $\Spec k$. The virtual tangent bundle
$T_f=[T_X]-[f^*T_Y]\in K_0(X)$ has rank $d=\dim X-\dim Y$.
Smooth purity for $p_Y$ identifies $\Sigma^\infty Y_+$ with
$p_{Y!}\Th_Y(T_Y)$ \cite[\S2.1.7]{DJK}. Since $f$ is proper,
the unit $B\to f_*f^*B$ and the identification $f_*=f_!$ give
a map $p_{Y!}B\to p_{X!}f^*B$ for every $B\in\SH(Y)$.
Taking $B=\Th_Y(T_Y)$, we define a transfer map $\tau_f$ as the composite
\begin{align*}
 \Sigma^\infty Y_+\simeq p_{Y!}\Th_Y(T_Y)
 &\longrightarrow p_{X!}\Th_X(f^*T_Y)\\
 &\simeq p_{X\#}\Th_X(-T_f)=\Th(-T_f).
\end{align*}
The last equivalence uses $f^*T_Y-T_X=-T_f$ and smooth purity on $X$.
For the same construction on smooth proper schemes, see
\cite[Remark~3.2]{NavarroGarmendiaDuality}.
Precomposition with $\tau_f$ commutes with maps of coefficient spectra.

\begin{lemma}\label{lem:transfer-comparison}
Let $f\colon X\to Y$ and $d$ be as above. For every $s\ge0$,
$x\in K_s(X)$ and $y\in\CH^j(X,s)$, one has
\begin{equation}\label{eq:pushforward-comparison}
 \beta^{-d}\tau_f^*\bigl(\ThomClass{\KGL}(-T_f)x\bigr)=f_!x,
 \qquad
 \tau_f^*\bigl(\ThomClass{\HZ}(-T_f)y\bigr)=f_*y.
\end{equation}
Here $f_!$ is the Thomason--Trobaugh direct image in algebraic $K$-theory,
induced by $Rf_*$ on perfect complexes, and $f_*$ is induced by proper
pushforward on Bloch's cycle complexes.
\end{lemma}

\begin{proof}
By smooth purity, the two composites give the direct image of
\cite[Definition~2.2]{NavarroNavarro}, expressed in Thom classes
as in \cite[Remark~3.2]{NavarroGarmendiaDuality}. The identifications
with the stated direct images follow from \cite[Corollary~2.12]{NavarroNavarro}
for $K$-theory and from the natural comparison
\cite[Corollary~8.12]{CDIntegral} for higher Chow groups.
\end{proof}

\needspace{3\baselineskip}
\begin{proposition}\label{prop:negative-weight}
\begin{enumerate}[label=\textup{(\alph*)}]
\item Let $Y$ be smooth of pure dimension $D$, let $Z\subset Y$ be closed
with $\dim Z<D-n$, and put $U=Y\setminus Z$. Then
$\CH^n(Y,s)\to\CH^n(U,s)$ is injective for every $s\ge0$.
\item Let $f\colon X\to Y$ and $d$ be as above. For every $n\ge0$ with
$n+d<0$, every $s\ge0$ and every $x\in K_s(X)$, one has
$\chi_{i,s}(f_!x)=0$ for every $0\le i\le n$.
In particular the last assertion of Theorem~\ref{thm:main} holds.
\end{enumerate}
\end{proposition}

\begin{proof}
\textup{(a)} Bloch's localization theorem for higher Chow groups \cite{BlochMoving},
as recalled in \cite[Properties~17.4(2)]{MVW},
identifies the kernel with the image of
$\CH_{D-n}(Z,s)$. Its cycle complex is zero, since in simplicial degree
$v$, its cycles would have dimension $D-n+v>\dim(Z\times\Delta_k^v)$.

\textup{(b)} Put $Z=f(X)$ and $D=\dim Y$. Then
$\dim Z\le\dim X=D+d<D-n\le D-i$.
Proper direct image commutes with restriction to $Y\setminus Z$,
whose inverse image in $X$ is empty. Thus $f_!x$ restricts to zero
there in every $K_s$.
Naturality of the Newton operations and part~(a)
give the result.
\end{proof}

\begin{theorem}\label{prop:basic-rr}
Let $k$ be a field of characteristic zero and let $f\colon X\to Y$ be a
projective morphism between smooth quasi-projective $k$-schemes of pure
dimension. Put $d=\dim X-\dim Y$ and $T_f=[T_X]-[f^*T_Y]$.
For every integer $n\ge0$ such that $m=n+d\ge0$, every integer
$M\ge\max(n,m)$, every $s\ge0$ and every $x\in K_s(X)$, one has
\begin{equation}\label{eq:basic-rr}
 \frac{T_M}{n!}\chi_{n,s}(f_!x)
 =f_*\bigl(T_M[\ch(x)\Td(T_f)]^{(s)}_m\bigr)
 \qquad\text{in }\CH^n(Y,s).
\end{equation}
\end{theorem}

\begin{proof}
Put $r=-d$ and
$z=\ThomClass{\kgl}(-T_f)x$.
The transfer gives
$\hat z=\tau_f^*z\in[\Sigma^{s,0}\Sigma^\infty Y_+,A_r]$,
mapping to $f_!x$ in $K_s(Y)$.

Suppose first that $d\ge0$, so $r\le0$.
Effective-cover adjunction identifies
\[
 [\Sigma^{s,0}\Sigma^\infty Y_+,A_0]
 \xrightarrow{\sim}
 [\Sigma^{s,0}\Sigma^\infty Y_+,A_r]
 \xrightarrow{\sim}K_s(Y).
\]
Thus $\hat z$ is the image of the effective lift of $f_!x$.
Apply $\psi^{(M)}_{r,r}$ and take position $m$.
By Lemma~\ref{lem:evaluation}, its value on $\hat z$ is
$(T_M/n!)\chi_{n,s}(f_!x)$, since $r+m=n$.
By Proposition~\ref{prop:thom}, its value on $z$ is
\[
 \ThomClass{\HZ}(-T_f)\,
 T_M[\ch(x)\Td(T_f)]^{(s)}_m.
\]
Naturality under $\tau_f^*$ and \eqref{eq:pushforward-comparison}
give the result.

If $d<0$, put $c=-d$. Instead apply
$\psi^{(M)}_{0,c}=\Phi_Mq_M\iota_{c,0}$ and take position $n$.
Its value on $\hat z$ is again the left side of \eqref{eq:basic-rr},
by Lemma~\ref{lem:evaluation}. Its value on $z$, after the Thom
isomorphism, is
$T_M[\ch(x)\Td(T_f)]^{(s)}_{n-c}$.
Since $n-c=m$, the same naturality argument concludes the proof.
\end{proof}

Taking $M=m$ proves Theorem~\ref{thm:main} when $d\ge0$.
\begin{theorem}
Let\label{prop:absolute} $m=n+d\ge0$ and $M\ge\max(m,n)$. Then
\[
 T_M[\ch(f_!x)\Td(T_Y)]^{(s)}_n
 =f_*\bigl(T_M[\ch(x)\Td(T_X)]^{(s)}_m\bigr).
\]
\end{theorem}
\begin{proof}
For $0\le i\le n$, put $M_i=\max(i,i+d)$.
Theorem~\ref{prop:basic-rr}, supplemented by the support vanishing
of Proposition~\ref{prop:negative-weight}(b) when $i+d<0$, gives
\[
 \frac{T_{M_i}}{i!}\chi_{i,s}(f_!x)
 =f_*\bigl(T_{M_i}[\ch(x)\Td(T_f)]^{(s)}_{i+d}\bigr),
\]
where a negative-weight term means zero. Since
$M_i+n-i=\max(n,m)\le M$, Lemma~\ref{lem:denominators} gives
$T_{M_i}T_{n-i}\mid T_M$. Multiply the weight-$i$ identity by
\[
 \frac{T_M}{T_{M_i}T_{n-i}}\,D_{n-i}(T_Y)
\]
and sum over $i$. The left side is the integral polynomial on the
left of the theorem. By the projection formula the right side is
the pushforward of the integral convolution
\[
 T_M\sum_{i=0}^n
 [\ch(x)\Td(T_f)]^{(s)}_{i+d}\,f^*\Td_{n-i}(T_Y).
\]
Putting $j=i+d$, this is the full weight-$m$ convolution if $d\le0$,
since negative-weight terms vanish. If $d>0$, its omitted terms
$0\le j<d$ have zero pushforward: each is the integral product
\[
 \frac{T_M}{T_jT_{m-j}}
 \bigl(T_j[\ch(x)\Td(T_f)]^{(s)}_j\bigr)
 f^*D_{m-j}(T_Y),
\]
and the first factor inside parentheses pushes forward to
$\CH^{j-d}(Y,s)=0$.
Thus one may use the full convolution in every case.
The identity $\Td(T_f)f^*\Td(T_Y)=\Td(T_X)$ identifies it with
$T_M[\ch(x)\Td(T_X)]^{(s)}_m$. This is an identity of integral
universal polynomials after clearing the displayed denominators,
so its evaluation does not require rational injectivity on $X$ or $Y$.
\end{proof}

\begin{unnumberedremark}[Relative and absolute forms]
For a fixed $f$, the families of identities in
Theorems~\ref{prop:basic-rr} and~\ref{prop:absolute}, with
$M\ge\max(n,n+d)$ and with the support vanishing for $n+d<0$,
are equivalent. For the reverse implication, use the weight-$i$ absolute identity
with multiplier $T_{M_i}$. If $i+d<0$, both sides are zero by
support vanishing. Multiply this identity by
$T_M D_{n-i}(-T_Y)/(T_{M_i}T_{n-i})$ and sum over $0\le i\le n$.
The preceding integral convolution argument applies with
$\Td(T_Y)^{-1}$, and cancels the absolute Todd factor.
\end{unnumberedremark}

\subsection{Negative relative dimension}\label{sec:negative}

Suppose $d=-c<0$ and $n\ge c$, so $m=n-c$. Recall from the proof
of Theorem~\ref{prop:basic-rr} that, for $x\in K_s(X)$, the class
\[
 \hat z=\tau_f^*\bigl(\ThomClass{\kgl}(-T_f)x\bigr)
 \in[\Sigma^{s,0}\Sigma^\infty Y_+,A_c]
\]
lifts $f_!x$ to the $c$-th effective cover $A_c=f_c\KGL$.
For $1\le r\le c$, the operation of \eqref{eq:psi} is the composite
\[
 \psi^{(n-r)}_{r,c}\colon
 A_c\xrightarrow{\ \iota_{c,r}\ }A_r\simeq\T^r A_0
 \xrightarrow{\ \T^r q_{n-r}\ }\T^r W_{n-r}
 \xrightarrow{\ \T^r\Phi_{n-r}\ }\bigoplus_{t=0}^{n-r}B_{r+t}.
\]
Its last component takes values in $B_n$. Naturality and
Proposition~\ref{prop:thom} give
\begin{equation}\label{eq:shifted-transfer}
 \psi^{(n-r)}_{r,c}(\hat z)_{n-r}
 =f_*\bigl(T_{n-r}[\ch(x)\Td(T_f)]^{(s)}_{n-c}\bigr).
\end{equation}
To express the left-hand side in terms of $\chi_{n,s}(f_!x)$,
we compare the two operations on $\Omega^\infty A_r$ for $r=1,2$.

\begin{lemma}\label{lem:low-covers}
Put $\mathcal K_r=\Omega^\infty A_r$. Let $BSL$ be the stable
classifying space of vector bundles with a trivialization of their
determinant. Then, in $\mathcal H_\bullet(k)$,
\[
 \mathcal K_1\simeq \BGL,\qquad
 \mathcal K_2\simeq BSL\simeq\operatorname{hofib}
       \bigl(\det\colon\BGL\longrightarrow B\Gm\bigr).
\]
For $r\in\{1,2\}$ and $R\in\{\Z,\Q\}$, the universal Chern classes give
a graded ring isomorphism
\[
 R[c_r,c_{r+1},\ldots]\xrightarrow{\ \sim\ }
 \bigoplus_{j\ge0}H^{2j}(\mathcal K_r,R(j)),\qquad \deg(c_i)=i.
\]
In particular, for each $j\ge0$ the coefficient map
$H^{2j}(\mathcal K_r,\Z(j))\to H^{2j}(\mathcal K_r,\Q(j))$
is injective.
\end{lemma}

\begin{proof}
Applying $\Omega^\infty$ to $A_1\to A_0\to\HZ$ gives the fiber of
the rank map $\Z\times BGL\to\Z$, hence $\mathcal K_1\simeq BGL$.
The next triangle gives the fiber of $BGL\to BG_m$, because
$\Omega^\infty\T\HZ=K(\Z(1),2)\simeq BG_m$.
This map rationalizes to $c_1$ and thus equals $c_1$ integrally,
since $H^2(BGL,\Z(1))=\Z c_1$. Its fiber is $BSL$; compare
\cite[Examples~4.4.2--4.4.3]{ABH}.

The $BGL$ calculation is in the proof of Lemma~\ref{lem:detection}.
For the integral cohomology of $BSL_e$, set
$G_N=\Gr(e,e+N)$, let $L_N=\det(\gamma_{e,N})$, and use
$P_N=L_N\setminus0$ as its approximation. The determinant torsor
models the $\mathbb A^1$-homotopy fiber of
$BGL_e\xrightarrow{\det}B\Gm$; in stable rank this is $BSL$
\cite[Example~4.4.3]{ABH}. In each fixed weight $j$, choose $N$
large enough that $\CH^a(G_N)=\Z[c_1,\ldots,c_e]_a$ for $a\le j$.
The localization sequence for $L_N$ and its zero section then gives
\[
 H^{2j}(P_N,\Z(j))=
 \CH^j(G_N)/c_1\CH^{j-1}(G_N),
\]
because $H^{2j-1}(G_N,\Z(j-1))=0$. In the preceding degree the
same sequence and injectivity of multiplication by $c_1$ in the
stable polynomial ring give
\[
 H^{2j-1}(P_N,\Z(j))=
 \bigl(\CH^{j-1}(G_N)/c_1\CH^{j-2}(G_N)\bigr)\otimes k^\times.
\]
These systems are eventually constant in each weight, so their
$\lim^1$ terms vanish. Passing first to $N$ and then to stable rank
gives $\Z[c_2,c_3,\ldots]$ on the diagonal. Repeating the calculation
with $\Q$ proves rational injectivity.
\end{proof}

\begin{theorem}\label{prop:improved}
Put $r=\min(c,2)$ and $L_{n,r}=\lcm(n!,T_{n-r})$. Then
\begin{equation}\label{eq:improved}
 \frac{L_{n,r}}{n!}\chi_{n,s}(f_!x)
 =f_*\bigl(L_{n,r}[\ch(x)\Td(T_f)]^{(s)}_{n-c}\bigr).
\end{equation}
\end{theorem}

\begin{proof}
Let $\epsilon_r\colon\Sigma^\infty\mathcal K_r\to A_r$ be the counit
and $u_r\colon\mathcal K_r\to\mathcal K$ the map induced by
$\iota_{r,0}$. The two integral classes
\[
 \frac{L_{n,r}}{T_{n-r}}
 \psi^{(n-r)}_{r,r}(\epsilon_r)_{n-r},
 \qquad \frac{L_{n,r}}{n!}u_r^*s_n
\]
have the same rationalization in
$\widetilde H^{2n}(\mathcal K_r,\Z(n))$.
Lemma~\ref{lem:low-covers} identifies them integrally.
Pull back along the adjoint of
$\iota_{c,r}\hat z\colon\Sigma^{s,0}\Sigma^\infty Y_+\to A_r$,
and use \eqref{eq:shifted-transfer}. The pulled-back Newton class is
$\chi_{n,s}(f_!x)$, including when $s>0$.
\end{proof}

For $c=1,2$ the multiplier in \eqref{eq:improved} is
$\lcm(n!,T_m)$. For every $c$ its prime factors are at most $n$, so
the formula recovers the coefficient range $\Z[1/n!]$ for $d<0$.

\begin{question}\label{question:uniform}
For a projective morphism of relative dimension $-c<0$, does
\eqref{eq:main} hold with $\lcm(n!,T_{n-c})$ for every $n\ge c$?
\end{question}

For $c\le2$ this is Theorem~\ref{thm:main}. For $c\ge3$, extending
the argument would require an integral identification of the
relevant character on $\Omega^\infty A_c$. The calculations for
$BGL$ and $BSL$ do not establish rational injectivity on these
higher covers; a rational equality there could have a torsion
difference.

For $3\le c<n\le6$, the two proposed coefficients coincide:
\[
 \lcm(n!,T_{n-2})=\lcm(n!,T_{n-c}).
\]
The first discrepancies are
\[
\begin{array}{c|c|c|c}
 n&c&\lcm(n!,T_{n-2})&\lcm(n!,T_{n-c})\\\hline
 7&3\le c\le6&10080&5040\\
 8&3\le c\le7&120960&40320.
\end{array}
\]
Thus the unresolved factor is $2$ in codimension $7$ and $3$ in
codimension $8$; the boundary case $m=0$ is already covered with
multiplier $n!$.

For the next cover, \cite[Example~4.4.3]{ABH} identifies
$\Omega^\infty A_3$, in characteristic zero, with the homotopy
fiber of $c_2:BSL\to K(\Z(2),4)$. This gives a concrete candidate
for extending the argument. The missing input is the required
integral-to-rational injectivity for the universal character classes
on this fiber;
the corresponding computation for $BSL$ does not establish it.

\subsection{Source weight zero}
\begin{proposition}\label{prop:weight-zero}
If $d=-c<0$ and $n=c$, then
\[
 \chi_{c,0}(f_!x)=c!f_*\rk(x),\qquad
 \chi_{c,s}(f_!x)=0\quad(s>0).
\]
Thus \eqref{eq:main} holds with $C_{c,-c}=c!$.
\end{proposition}

\begin{proof}
Put $D=\dim Y$, $t=D-c=\dim X$ and $Z=f(X)_{\mathrm{red}}$.
Since $f_!x$ restricts to zero on $Y\setminus Z$, localization
places $\chi_{c,s}(f_!x)$ in the image of
$\CH_t(Z,s)\to\CH^c(Y,s)$. As $\dim Z\le t$, the only generators
of $z_t(Z,v)$ are $Z_j\times\Delta_k^v$, where $Z_j$ runs over the
$t$-dimensional irreducible components of $Z$. They form a constant
simplicial abelian group, so $\CH_t(Z,s)=0$ for $s>0$.

Now let $s=0$. By Proposition~\ref{prop:negative-weight}(a) we may
discard from $Y$ any closed subset of dimension less than $t$.
We may thus assume that $Z$ is the disjoint union of the $Z_j$.
Over their generic points the map $g\colon X\to Z$ is finite, since
$\dim X=t$. By properness and generic flatness, after a further such
restriction $g$ is finite flat. Shrinking each $Z_j$ once more, we may
trivialize the vector bundles representing $g_!x$ there. On the
resulting open $U\subset Y$, write $Z_j^\circ=Z_j\cap U$ and let
$r_j$ be the generic rank of $g_!x$ on $Z_j^\circ$. Then
\[
 (f_!x)|_U=\sum_j r_j[\mathcal O_{Z_j^\circ}],\qquad
 (f_*\rk(x))|_U=\sum_j r_j[Z_j^\circ],
\]
where the second equality is the definition of proper pushforward
of cycles at the generic points. By the calculation following
Lemma~\ref{lem:slice-normalization},
$\chi_{c,0}([\mathcal O_{Z_j^\circ}])=c![Z_j^\circ]$. Additivity
gives the formula on $U$, and injectivity of restriction gives it on $Y$.
\end{proof}

\begin{proof}[Proof of Theorem~\ref{thm:main}]
For $m\ge0$, the formula follows from Theorems~\ref{prop:basic-rr}
and~\ref{prop:improved} and Proposition~\ref{prop:weight-zero}.
For $m<0$, the vanishing is Proposition~\ref{prop:negative-weight}(b).
\end{proof}

\section{The characteristic of the ground field}\label{sec:characteristic}

Let $k$ be a field of exponential characteristic $e$. We first prove
the Riemann--Roch formula over a perfect field with coefficients in
$\Z[1/e]$ and the same multipliers $C_{n,d}$ as in characteristic zero.
We then pass to an arbitrary field by descent from its perfect closure.

\begin{remark}\label{rem:integral-perfect-cases}
Over a perfect field the cases $m\le0$ of Theorem~\ref{thm:main} hold
integrally. For $m<0$, the localization argument of
Proposition~\ref{prop:negative-weight}(b) works over any field.
For $m=0$ and $n>0$, the proof of Proposition~\ref{prop:weight-zero}
uses only localization, generic ranks and the integral formula for
$\chi_{n,0}([\mathcal O_Z])$. When $m=n=0$, the identity in degree
$s=0$ is $\rk(f_!x)=f_*\rk(x)$ at the generic points of $Y$, and
$\CH^0(Y,s)=0$ for $s>0$.
\end{remark}

\begin{theorem}\label{thm:positive-characteristic}
Let $k$ be a perfect field of exponential characteristic $e$, and put
$\Delta=\Z[1/e]$. The integral realization theorem holds in
$\SH(k)[1/e]$, with $B_q=\T^q\mathrm H \Delta$.
The image of top restriction
$[W_q[1/e],\T^q\mathrm H \Delta]\to \Delta$ is $T_q\Delta$;
for $q\ge1$, the additive order of $\kappa_q[1/e]$ is the part of
$T_q$ prime to $e$. The Riemann--Roch identities of
Theorem~\ref{thm:main} hold in $\CH^n(Y,s)\otimes \Delta$ for
$x\in K_s(X)\otimes \Delta$, with the same integers $C_{n,d}$.
\end{theorem}

\begin{proof}
The slice identifications for $\KGL$ over perfect fields follow from
\cite[Theorems~6.4.2 and~9.0.3]{Levine}; after inverting $e$ they also
follow from the Landweber slice calculation
\cite[Theorem~8.5]{Hoyois}. The zero slice of the sphere with these
coefficients is $\mathrm H \Delta$; see
\cite[Theorem~10.5.1]{Levine} or \cite[Theorem~8.7]{Hoyois}.
The stable Adams operations are defined over regular bases after
inverting their index \cite[Definition~5.3.2 and the following
paragraph]{Riou}; the rational character is defined over perfect
fields \cite[\S6.2.3]{Riou}.

In the primewise argument take only primes
$\ell\ne\operatorname{char}k$. In this range the dual Steenrod algebra
has the same Milnor description \cite[Theorem~1.1 and~\S5.1]{HKO},
\cite[Theorem~5.17]{Hoyois}, so the Bockstein calculation of
Proposition~\ref{prop:bockstein-basis} applies. Since $e$ is invertible
in $\Z_{(\ell)}$, the spectrum $\mathrm H\Z_{(\ell)}$ is cellular
\cite[Proposition~8.1]{Hoyois}. The proof of Proposition~\ref{cc:prop:modules}
therefore gives the module equivalence in this range. As in
Corollary~\ref{cc:cor:odd}, the odd obstruction groups are killed by
$\ell$, and the argument of Proposition~\ref{cc:prop:prime} gives the exponent
$\lfloor q/(\ell-1)\rfloor$.

Set $S=\{\ell:\ell\mid T_q,\ \ell\nmid e\}$ and
$D=\Delta[1/T_q]$. The analogue of Lemma~\ref{cc:lem:arithmetic}
over $\Delta$ follows from the exact sequence
\[
 0\longrightarrow \Delta\longrightarrow
 D\oplus\bigoplus_{\ell\in S}\Z_{(\ell)}
 \longrightarrow\bigoplus_{\ell\in S}\Q\longrightarrow0.
\]
The projector of Proposition~\ref{cc:prop:D} over $D$ and the
$\ell$-local projections glue as in
Subsection~\ref{cc:sec:glue}. If $S$ is empty, $D=\Delta$ and the projector
alone suffices. The exact-denominator argument of
Theorem~\ref{thm:exact-denominator} works over $\Delta$ and gives the
ideal $T_q\Delta$, and hence the asserted order of the attaching map.

For the geometric comparison, evaluate localized operations on
the integral counits of $\Z\times BGL$, $BGL$ and $BSL$, followed
by localization. The first two effective-cover identifications
remain valid by the integral slice calculation over perfect fields.
The Grassmannian and determinant-line computations in
Lemmas~\ref{lem:detection} and~\ref{lem:low-covers} hold with
$\Delta$-coefficients: diagonal groups on each rank component are free,
and the preceding inverse systems are surjective. Products over
the rank components inject componentwise into the corresponding
rational groups. The evaluations of Lemma~\ref{lem:evaluation} and
Proposition~\ref{prop:thom}, followed by naturality under the proper
transfer, therefore prove the identities for integral
$K$-classes with values in $\CH^n(Y,s)\otimes \Delta$. Additivity gives
the assertion for $K_s(X)\otimes \Delta$. The cases $m\le0$ also follow
from Remark~\ref{rem:integral-perfect-cases}.
\end{proof}

\begin{lemma}\label{lem:perfection-injective}
Let $k$ be any field of exponential characteristic $e$, let
$k^{\mathrm{perf}}$ be its perfect closure, and let $Y$ be an
equidimensional scheme of finite type over $k$. Then
\[
 \CH^n(Y,s)\otimes\Z[1/e]\longrightarrow
 \CH^n(Y_{k^{\mathrm{perf}}},s)\otimes\Z[1/e]
\]
is injective.
\end{lemma}

\begin{proof}
There is nothing to prove in characteristic zero. Otherwise write
$k^{\mathrm{perf}}=\varinjlim_L L$, where $L/k$ ranges over finite
purely inseparable extensions. Flat pullback gives an isomorphism
of Bloch cycle complexes
\[
 \varinjlim_L z^n(Y_L,\bullet)
       \xrightarrow{\ \sim\ }
 z^n(Y_{k^{\mathrm{perf}}},\bullet).
\]
Indeed, the ideal of an integral closed subscheme in a simplex
descends to some finite $L$. Its descended subscheme is integral
by faithful flatness, and its pullback is the original reduced
subscheme, so the cycle pulls back with multiplicity one.
The face conditions are preserved and reflected by field
extension, since the relevant dimensions are unchanged.
Conversely, flat pullback is injective on cycle groups: distinct
integral supports have disjoint sets of generic points after
base change, and every pullback has positive multiplicities.
This proves the displayed isomorphism degree by degree, compatibly
with the differentials. Exactness of filtered colimits gives the
same continuity statement for higher Chow groups.

For $\pi_L:Y_L\to Y$, proper pushforward and flat pullback satisfy
\[
 (\pi_L)_*\pi_L^*=[L:k]\,\mathrm{id}
\]
already on cycle complexes; see \cite[\S4, equation~(4.1)]{LevineCycles}.
Since $[L:k]$ is a power of $e$, every $\pi_L^*$ is injective
after inverting $e$. Continuity now proves the assertion.
\end{proof}

\begin{corollary}\label{cor:arbitrary-field}
Over an arbitrary field of exponential characteristic $e$,
Theorem~\ref{thm:main} holds with $\Z[1/e]$-coefficients and the
same integers $C_{n,d}$.
\end{corollary}

\begin{proof}
Base change to $k^{\mathrm{perf}}$. Smoothness, projectivity and
dimensions are preserved, and
$T_{f_{k^{\mathrm{perf}}}}$ is the pullback of $T_f$.
Flat base change identifies the pullback of $Rf_*$ with the
direct image after base change. Proper pushforward of higher
cycles commutes with flat pullback, and the Chern operations are
natural under field extension. Thus the difference between the
two sides of each asserted identity pulls back to zero by
Theorem~\ref{thm:positive-characteristic} and
Remark~\ref{rem:integral-perfect-cases}.
Lemma~\ref{lem:perfection-injective} makes that difference zero
over $k$. This proves the formulas for integral $K$-classes;
additivity extends them to $K_s(X)\otimes\Z[1/e]$.
\end{proof}

The field hypothesis enters the vanishing
$H^a(k,\Delta(b))=0$ for $a>b$ and for $b<0$, used in Lemma~\ref{cc:lem:rational}, in the computation of odd
operations in Corollary~\ref{cc:cor:odd}, and in the universal cohomology
calculations. Quasi-projectivity supplies Jouanolou's affine
bundle torsor and the usual vector-bundle presentations. In this
setting every proper morphism is projective, since an ample line
bundle on the source is relatively ample. The geometric arguments
use no resolution of singularities.

\begin{remark}\label{rem:motivic-cycle-pushforward}
Let $f\colon X\to Y$ be a projective morphism between smooth
quasi-projective $k$-schemes as in Theorem~\ref{thm:main}.
The motivic Gysin map of Section~\ref{sec:arithmetic-rr} is defined
integrally, and for smooth $U$ over $k$ one has
$H^{2j-s}(U,\Z(j))\simeq\CH^j(U,s)$
\cite[Theorem~1.1(8)]{BEM}. The comparison of the Gysin map defined
by duality with the oriented projective Gysin map is given by
\cite[Corollary~2.6]{NavarroNavarro} and \cite[\S4.4.3]{DJK}.
For $s>0$, the identification of this map with proper pushforward on
Bloch's higher cycles used here has coefficients in $\Z[1/e]$
\cite[Theorem~8.11 and Corollary~8.12]{CDIntegral}.
For $s=0$ and $k$ perfect, the Gysin map agrees integrally with the
usual pushforward on Chow groups \cite[Example~3.2.9(1)]{DegliseOrientation}.
Thus, for $m>0$, Theorem~\ref{thm:arithmetic-rr} and
Corollary~\ref{cor:arithmetic-shift} give the $K_0$ formula with
multiplier $C_{n,d}$ in $\CH^n(Y)$ over every perfect field.
\end{remark}

\section{Integral Riemann--Roch over arbitrary schemes}\label{sec:arithmetic-rr}

Let $S=\Spec\Z$. We use Spitzweck's oriented spectrum $\HZ_S$ and,
for a quasi-compact quasi-separated scheme $g\colon U\to S$, put
$\HZ_U=g^*\HZ_S$. Over $S$ the identification
$\HZ_S\simeq s_0\KGL_S$ follows from Bachmann's description of the
effective cover and its Bott cofiber
\cite[Theorem~1.1(1), Example~1.3 and Lemma~2.1]{BachmannDedekind}.
Bachmann, Elmanto and Morrow prove that the comparison maps
\[
 g^*\kgl_S\xrightarrow{\sim} f_0\KGL_U,\qquad
 g^*s_0\KGL_S\xrightarrow{\sim}s_0\KGL_U
\]
are equivalences \cite[Remark~4.14 and Corollaries~9.6--9.7]{BEM}.
Thus the cohomology used here is their $\mathbb A^1$-invariant
motivic cohomology. Our conventions are
\[
 H^i(U,\Z(j))=[\one_U,\Sigma^{i,j}\HZ_U]_{\SH(U)},\qquad
 KH_s(U)=[\Sigma^{s,0}\one_U,\KGL_U]_{\SH(U)}.
\]
Here $KH$ is homotopy invariant $K$-theory. For a regular noetherian
scheme of finite dimension the comparison $K_s(U)\to KH_s(U)$ is an
isomorphism.

\begin{unnumberedremark}
If $U$ is smooth over a field, then
$H^{2j-s}(U,\Z(j))\simeq\CH^j(U,s)$ \cite[Theorem~1.1(8)]{BEM}.
On a scheme smooth over a mixed-characteristic Dedekind domain the
represented theory is the Zariski hypercohomology of the Bloch--Levine
cycle complex \cite[Corollary~7.19]{Spitzweck}, \cite[Remark~4.14]{BEM}.
For smooth arithmetic schemes this gives $H^{2j}(U,\Z(j))=\CH^j(U)$
and $H^t(U,\Z(j))=0$ for $t>2j$ \cite[Theorem~3.2]{KahnCycle}.
In particular, for $j\ge0$,
\begin{equation}\label{eq:arithmetic-base-coefficients}
 H^{2j}(S,\Z(j))=
 \begin{cases}\Z&j=0,\\0&j>0,\end{cases}
 \qquad H^{2j+1}(S,\Z(j))=0.
\end{equation}
Indeed $\operatorname{Pic}(S)=0$ and $\dim S=1$. These vanishings
on the base will determine a universal operation before it is
evaluated on a possibly singular scheme.
\end{unnumberedremark}

\begin{theorem}\label{thm:arithmetic-rr}
Let $f\colon X\to Y$ be proper between quasi-compact quasi-separated
schemes. Suppose it factors as $X\xrightarrow{i}P\xrightarrow{p}Y$,
where $i$ is a finitely presented regular closed immersion and $p$ is
smooth, separated and of finite presentation. Assume that
$T_f=i^*T_{P/Y}-N_i$ has constant rank $d$. Let $f_!$ and $f_*$
denote the oriented Gysin direct images in $KH$ and in motivic
cohomology. For $n,s\ge0$, put $m=n+d$ and $M=\max(n,m)$.
For every $x\in KH_s(X)$ one has
\begin{equation}\label{eq:arithmetic-rr}
 \frac{T_M}{n!}\,\chi_{n,s}(f_!x)
 =f_*\!\left(T_M[\ch(x)\Td(T_f)]^{(s)}_m\right)
 \quad\text{in }H^{2n-s}(Y,\Z(n)).
\end{equation}
For $m<0$ the right side is zero. For $m\ge0$ the class inside
$f_*$ is the integral polynomial
\[
 \sum_{i=0}^m\frac{T_M}{i!T_{m-i}}\,
       \chi_{i,s}(x)D_{m-i}(T_f)\ \in H^{2m-s}(X,\Z(m)).
\]
Here $\chi_{i,s}$ is the universal Newton operation on $KH$ with the
normalization of \eqref{eq:newton}. In particular the identity applies
to the image of any class under $K_s(X)\to KH_s(X)$.
\end{theorem}

The identity with motivic Gysin maps is therefore integral over
fields of positive characteristic. For the comparison with proper
pushforward on Bloch's higher cycles, Section~\ref{sec:characteristic}
uses coefficients in which the exponential characteristic is invertible.

\subsection{The stable character over the integers}

Write $A_{i,U}=f_i\KGL_U$ and $B_{i,U}=s_i\KGL_U$. Bott periodicity
identifies them with $\T^i\kgl_U$ and $\T^i\HZ_U$, respectively.
The tower maps are the Bott translates of $\T\kgl_U\to\kgl_U$.
Consequently the covers, the slices and the maps between them are
pulled back from $S$. Over $S$ we omit the subscript and set
$W_q=\cofib(A_{q+1}\to A_0)$. We shall construct the operation over
this single base.

\begin{lemma}\label{lem:arithmetic-character}
For $q\ge0$ there is a map $\theta_q\colon W_q\to B_q$ such that
$\theta_q\iota_q=T_q\id_{B_q}$ and $(\theta_q)_\Q=T_q\ch_q$.
For $q=0$ one can take the identity.
\end{lemma}
\begin{proof}
For $q=0$ take the identity of $W_0=B_0$. Suppose $q>0$.
We verify the two inputs to the extension argument of
Subsections~\ref{cc:sec:operations}--\ref{cc:sec:glue}. Fix a prime $p$.
Write $R=\HZ_{(p)}$, $H=\HF_p$ and $H_{**}=\pi_{**}H$ in $\SH(S)$.
By \cite[Theorems~7.1--7.2]{BachmannHopkinsCharacteristic}, the
Milnor monomials give a free right $H_{**}$-basis of
$\pi_{**}(H\wedge H)$, and $\HZ_S\simeq\mathrm{MGL}_S/(t_1,t_2,\ldots)$
is cellular. Both statements apply over $S$, including the fiber of
characteristic $p$.

Let $d_L,d_R$ be the Bocksteins on the two factors of $H\wedge H$.
For the left and right images of the orientation classes $u,v$,
Spitzweck's coaction formulas are
\[
 u_L=u_R+\sum_{i\ge0}\tau_i v_R^{p^i},\qquad
 v_L=\sum_{i\ge0}\xi_i v_R^{p^i},\qquad \xi_0=1,
\]
with $\beta u=v$ and $\beta v=0$
\cite[Lemma~10.6 and Corollary~10.25]{Spitzweck}. The basis map is
the one in \cite[Theorem~10.26]{Spitzweck}, whose equivalence at
the residue characteristic is supplied by Bachmann--Hopkins.
Applying $d_L$ and $d_R$ separately and comparing coefficients gives
\[
 \begin{gathered}
 d_L\tau_i=\xi_i\ (i\ge0),\qquad
 d_R\tau_0=-1,\qquad d_R\tau_i=0\ (i>0),\\
 d_L\xi_i=d_R\xi_i=0 .
 \end{gathered}
\]
Reduction $R\wedge H\to H\wedge H$ identifies
$\pi_{**}(R\wedge H)$ with $\ker d_L$. Indeed multiplication by $p$
on $R\wedge H$ is null, so reduction is injective, and the connecting
map followed by reduction is $d_L$. To compute this kernel, write an
element uniquely as $v_0+\tau_0 w_0$ on the ordered basis, with
$v_0,w_0$ free of $\tau_0$. On these monomials put $\delta=d_L$.
Since $d_L$ kills right coefficients,
\[
 d_L(v_0+\tau_0 w_0)=\delta v_0+w_0-\tau_0\delta w_0.
\]
Thus the elements $b'=b-\tau_0\delta b$, for basis monomials $b$
free of $\tau_0$, form a right coefficient basis of the kernel.
On these generators $d_R b'=\delta b=(\delta b)'$.
On coefficients $d_R$ is the coefficient Bockstein. After exchanging
factors, the reduced part is therefore the complex of
Lemma~\ref{cc:lem:pairs}, with differential $\tau_i\mapsto\xi_i$
for $i>0$, tensored with the coefficient complex.
This uses the ordered basis also at $p=2$.

To pass from cooperations to operations, let $C=\fib(R\wedge R\to R)$.
The rational multiplication is an equivalence since
$(\HZ_S)_\Q$ is Beilinson motivic cohomology
\cite[Theorem~7.14]{Spitzweck}. The preceding reduced complex is
acyclic, and the long exact sequence of $C\xrightarrow{p}C\to C/p$
therefore gives $p\pi_{**}C=0$ by Lemma~\ref{cc:lem:coop}.
Cellularity now allows Proposition~\ref{cc:prop:modules} to realize
the Koszul pairs as a splitting of $C$ into suspensions of $H$ as
$R$-modules. Applying $[-,\Sigma^{2h+1,h}R]_R$ to $R\wedge R=R\oplus C$,
the first summand contributes $H^{2h+1}(S,\Zp(h))=0$ by
\eqref{eq:arithmetic-base-coefficients}. The other summands contribute
groups killed by $p$. Free-module adjunction gives
\begin{equation}\label{eq:arithmetic-odd}
 p\,[R,\Sigma^{2h+1,h}R]_{\SH(S)}=0\qquad(h\ge0).
\end{equation}

For the second input, $\Psi^u$ on $\KGL[1/u]$ preserves the unit and
takes $\beta$ to $u\beta$
\cite[\S3.3.2, Lemma~3.11 and Proposition~3.12]{BachmannHopkinsAdams}.
The unit induces $s_0\one\simeq B_0$ \cite[Corollary~9.7]{BEM},
so $\Psi^u$ acts on $B_i[1/u]$ by $u^i$. Choose $u$ prime to $p$
generating $\Fp^\times$ modulo $p$. The Adams summand containing
$B_{q,(p)}$ has successive weights $q-t(p-1)$ for
$0\le t\le h_p=\lfloor q/(p-1)\rfloor$. At step $t>0$ the
obstruction to extending the preceding map lies in
\[
 [\Sigma^{-1}B_{q-t(p-1),(p)},B_{q,(p)}]
   =[R,\Sigma^{2t(p-1)+1,t(p-1)}R].
\]
By \eqref{eq:arithmetic-odd}, multiplying that map by $p$ kills the
obstruction. Induction as in Proposition~\ref{cc:prop:prime} extends
$p^{h_p}\id$ over $W_{q,(p)}$. Over $\Z[1/T_q]$ the finite
difference in Proposition~\ref{cc:prop:D} extends the top identity.

For $q>0$, scale the map over $\Z[1/T_q]$ by $T_q$ and the map at
$p\mid T_q$ by $T_q/p^{h_p}$. They agree on the top slice. For $i<q$
we have
\[
 [B_i,B_{q,\Q}]=H^{2(q-i)}(S,\Q(q-i))=0
\]
by rational idempotence and \eqref{eq:arithmetic-base-coefficients}.
Induction on the finite tower gives $[W_{q-1},B_{q,\Q}]=0$.
Hence a rational map from $W_q$ to $B_{q,\Q}$ is determined by its
top restriction. The rational Chern projection has top value one
by the rank and Bott normalization
\cite[Theorem~4.48, Corollary~4.52 and Remark~4.53]{BEM}.
The projection factors through $W_q$ since $B_{q,\Q}$ is right
orthogonal to $(q+1)$-effective objects. Thus the scaled maps all
rationalize to $T_q\ch_q$.
Applying $[W_q,-]$ to the coefficient pullback of
Lemma~\ref{cc:lem:arithmetic} gives an integral lift. Its top
restriction is determined by its rational value, since slice
orthogonality and $s_0\one\simeq\HZ_S$ give
$[B_q,B_q]=[\one_S,\HZ_S]=\Z$. Hence it equals $T_q\id$.
\end{proof}

\begin{unnumberedremark}
There is a shorter construction with a larger multiplier. For $q>0$,
$P_u(t)=\prod_{i=0}^{q-1}(t-u^i)$ annihilates $\Psi^u$ on $W_{q-1}$
after inverting $u$, by Lemma~\ref{cc:lem:polynomial}. Naturality
of the attaching map gives
\[
 P_u(u^q)\kappa_q=0,\qquad
 P_u(u^q)=u^{q(q-1)/2}\prod_{j=1}^q(u^j-1).
\]
Taking $u=2$ over $\Z[1/2]$ and $u=3$ over $\Z_{(2)}$ shows that
$N_q=\prod_{j=1}^q(2^j-1)2^{v_2(3^j-1)}$ kills the attaching map
in both localizations. The lifts of $N_q\id$ glue by the same rational
uniqueness. Already $N_2=48$, whereas $T_2=12$.
\end{unnumberedremark}

\subsection{Evaluation on Thom classes and the Gysin map}

Fix $M\ge0$. Compose $A_{0,S}\to W_M$ with $\theta_M$ and pull the
composite back to $U$. For integers $a\le b$ denote the resulting
operation by
\[
 F^{(M)}_{a,b,U}\colon A_{b,U}\xrightarrow{\iota_{b,a}} A_{a,U}
 \simeq\T^a A_{0,U}\longrightarrow \T^{a+M}\HZ_U .
\]
Its rationalization is the restriction of $T_M\ch_{a+M}$ by Bott
periodicity. We need its integral value both on a Thom class over $X$
and on an untwisted class over $Y$.

\begin{lemma}\label{lem:arithmetic-evaluation}
Let $\xi=V-W$ be a difference of vector bundles of constant virtual
rank $r$ on $U$, and let $x\in KH_s(U)$, $s\ge0$. Using Bott
periodicity, the class $\ThomClass{\KGL}(\xi)x$ has a unique lift
\[
 z\in[\Sigma^{s,0}\Th_U(\xi),A_{r,U}].
\]
For $a\le r$ put $j=a+M-r$. Under the Thom isomorphism
\[
 [\Sigma^{s,0}\Th_U(\xi),B_{a+M,U}]\simeq H^{2j-s}(U,\Z(j)),
\]
the class $F^{(M)}_{a,r,U}z$ corresponds to
$T_M[\ch(x)\Td(-\xi)]^{(s)}_j$ for $j\ge0$. It is zero for $j<0$.
\end{lemma}
\begin{proof}
Locally on a simultaneous trivializing cover,
$\Th_U(\xi)\simeq\T^r\one_U$. Descent for this cover expresses the
Thom object by its restrictions to the intersections and the smooth
pushforwards along the open immersions. These pushforwards preserve
$r$-effectivity, as does simplicial suspension. Effective-cover
adjunction therefore gives an equivalence of mapping spaces into
$A_{r,U}$ and into $\KGL_U$ from this source. It provides the lift
and its pullback naturality. If $j<0$, then $a+M<r$ and
$[A_{r,U},B_{a+M,U}]=0$ by slice orthogonality. We may assume
$0\le j\le M$, where the upper bound follows from $a\le r$.

First let $s=0$. Work on each open and closed piece on which
$e=\rk V$ and $h=\rk W$ are fixed. Applying $F^{(M)}_{a,r}$ and
then the inverse $\HZ$-Thom isomorphism gives an operation in
$(V,W,x)$, pointed at $x=0$. The Thom construction is functorial in
bundle isomorphisms and pullbacks \cite[\S16.2]{BachmannHoyois}.
Together with the mapping-space adjunction above, this gives an
operation of spaces. The polynomial \eqref{eq:integral-polynomial}
with $\xi$ replaced by $-\xi$ gives a second such operation.

We describe their common domain of comparison. On the site of
quasi-compact quasi-separated schemes let $L$ denote localization
for cdh descent and $\mathbb A^1$-invariance. Write $\mathrm{Vect}_e$
for the groupoid of rank-$e$ bundles and $\cK=\Omega^\infty KH$.
The tautological maps give, after applying $L$, an equivalence
\[
 \Gr(e,\infty)\times\Gr(h,\infty)\times(\Z\times\Gr)
 \ \longrightarrow\ \mathrm{Vect}_e\times\mathrm{Vect}_h\times\cK .
\]
Indeed $B\mathrm{GL}_e$ presents $\mathrm{Vect}_e$ Zariski locally,
and the infinite frame space gives its Grassmannian model.
Group completion of the space of bundles gives the $K$-theory
space Zariski locally, and its cdh sheafification is
$\Omega^\infty KH$ \cite[Remark~16.11, Lemma~16.12 and the discussion
following it]{BachmannHoyois}. To pass the Grassmannian models from
smooth $S$-schemes to this site, left Kan extend their presentations.
The schemes occurring in the Grassmannians and the bar constructions
for $\mathrm{GL}_e$ are all smooth over $S$. This extension preserves
products and carries the generating Nisnevich squares and affine-line
equivalences to local equivalences on the larger site. The group
completion comparison there follows on affines, where exact
sequences of projective modules split, and then by Zariski descent.
The cohomology target is already local for $L$
\cite[Proposition~3.30 and Theorem~1.1]{BEM}.

It remains to compare the universal classes. Let $G$ be a finite
Grassmannian product in the displayed model, on a fixed rank component
of $\Z\times\Gr$. For a tuple $\lambda$ of Schur indices let
$|\lambda|$ be its total codimension. The Grassmann bundle formula
\cite[Propositions~6.1--6.2]{Landweber} gives
\[
 H^{t}(G,\Z(j))\simeq
 \bigoplus_{\lambda}H^{t-2|\lambda|}(S,\Z(j-|\lambda|)).
\]
For $t=2j$, negative weights vanish by slice orthogonality, and
\eqref{eq:arithmetic-base-coefficients} leaves a copy of $\Z$ for
each $|\lambda|=j$. Restrictions to smaller stages are surjective
also for $t=2j-1$, since each Schur class at the smaller stage lifts.
The Milnor $\lim^1$ therefore vanishes on a cofinal sequence of
stages. It follows, component by component, that rationalization
detects the integral diagonal classes of the universal space. The
same holds relative to $x=0$, since the inclusion at $x=0$ is split
by projection onto the two bundle factors.

Rationally the classes agree by \eqref{eq:rational-thom}. For a line
of first Chern class $v$ the ratio of the two Thom classes is
$(1-\exp(-v))/v$. Multiplicativity and the splitting principle
give $\Td(\xi)^{-1}=\Td(-\xi)$ for $V-W$. The coefficient bound
$j\le M$ makes the polynomial integral by \eqref{eq:divisibility}.
The preceding detection proves equality as operations pointed in $x$,
so the equality can be evaluated on every $U$.

For $s>0$ apply $\pi_s$ at $x=0$ in the $KH$ variable, keeping
$V,W$ fixed. Products of two positive Chern classes of $x$ pull
back to zero on the simplicial suspension because its reduced
diagonal is null. The rank term also vanishes. Hence the degree-$i$
Newton polynomial contributes $(-1)^{i-1}i\,c_{i,s}(x)$ for $i>0$,
and the Todd coefficients are carried along unchanged. This gives
the formula with the normalization \eqref{eq:newton}.
\end{proof}

\begin{proof}[Proof of Theorem~\ref{thm:arithmetic-rr}]
Put $r=-d$, $a=n-M$ and $\xi=-T_f$. Then $a\le\min(0,r)$,
$a+M=n$ and $a+M-r=m$. The fundamental class of $f$ gives, for
every $E\in\SH(Y)$, the purity transformation
$\mathfrak p_f(E)\colon\Th_X(T_f)\otimes f^*E\to f^!E$
\cite[Theorem~3.3.2 and \S4.3.1]{DJK}. Define
\[
 G_f^E\colon[\Sigma^{s,0}\Th_X(\xi),f^*E]
            \longrightarrow[\Sigma^{s,0}\one_Y,E]
\]
by the composite
\begin{align*}
 [\Sigma^{s,0}\Th_X(\xi),f^*E]
 &\simeq [\Sigma^{s,0}\one_X,\Th_X(T_f)\otimes f^*E]
   \xrightarrow{\mathfrak p_f(E)}[\Sigma^{s,0}\one_X,f^!E]\\
 &\simeq [\Sigma^{s,0} f_!\one_X,E]
   \longrightarrow[\Sigma^{s,0}\one_Y,E].
\end{align*}
The second equivalence is adjunction. The last arrow precomposes with
the suspension of the unit $\one_Y\to f_*\one_X\simeq f_!\one_X$,
using properness. The six operations are taken relative to $Y$, so
the finite presentation hypothesis is imposed on $f$.
For $E=\KGL_Y$ and $E=B_{n,Y}$, the Thom isomorphisms identify
$G_f^E$ with the two oriented direct images \cite[\S4.3.3(iii)]{DJK}.

Write $D_X=\Sigma^{s,0}\Th_X(\xi)$ and
$D_Y=\Sigma^{s,0}\one_Y$. Naturality in $E$, together with
$f^*F^{(M)}_{a,r,Y}=F^{(M)}_{a,r,X}$, gives the square
\[
\begin{tikzcd}[column sep=large,row sep=large]
 {[D_X,A_{r,X}]} \arrow[r,"F^{(M)}_{a,r,X}"] \arrow[d,"G_f^{A_{r,Y}}"']
 & {[D_X,B_{n,X}]} \arrow[d,"G_f^{B_{n,Y}}"] \\
 {[D_Y,A_{r,Y}]} \arrow[r,"F^{(M)}_{a,r,Y}"']
 & {[D_Y,B_{n,Y}]}
\end{tikzcd}
\]
Take $z$ from Lemma~\ref{lem:arithmetic-evaluation} and write
$\hat z=G_f^{A_{r,Y}}z$. Naturality for $A_{r,Y}\to\KGL_Y$ shows
that $\hat z$ maps to $f_!x$. Since $D_Y$ is effective and $a\le0$,
$[D_Y,A_{a,Y}]\to[D_Y,\KGL_Y]$ is an isomorphism. Thus
$\iota_{r,a}\hat z$ and the effective lift of $f_!x$ followed by
$\iota_{0,a}$ coincide there. The lower horizontal map evaluates to
$(T_M/n!)\chi_{n,s}(f_!x)$ by the lemma with $\xi=0$.
The upper horizontal map evaluates, after Thom isomorphism, to
$T_M[\ch(x)\Td(T_f)]^{(s)}_m$, or to zero if $m<0$.
Commutativity proves \eqref{eq:arithmetic-rr}.
\end{proof}

\begin{samepage}
\begin{unnumberedremark}
If the factorization in Theorem~\ref{thm:arithmetic-rr} can be chosen
with $P=\mathbf P_Y(E)$ for a vector bundle $E$ and $p$ its projection,
then $f_!$ on $KH_s$ agrees with the map induced by $Rf_*$ on perfect
complexes. The uniqueness argument of \cite[Theorem~2.32]{Navarro}
applies to this factorization. It is routine to check that both direct
images satisfy functoriality, normalization for Cartier divisors, the
excess formula for blow-ups and the projection formula.
For a proper morphism between smooth schemes over a regular noetherian
base of finite dimension, the comparison in $K$-theory is proved in
\cite[\S2.1, Corollary~2.12]{NavarroNavarro}.
\end{unnumberedremark}
\end{samepage}

\subsection{Negative relative dimension and smooth arithmetic schemes}

\begin{corollary}\label{cor:arithmetic-shift}
Suppose $d=-c<0$. For every $n\ge0$ put $b=\min(c,n,2)$ and
$L=\lcm(n!,T_{n-b})$. Then \eqref{eq:arithmetic-rr} holds with $L$
in place of $T_M$. For $m>0$ this is $C_{n,d}$ of \eqref{eq:multiplier}.
\end{corollary}
\begin{proof}
On the same site consider the pointed spaces
\[
 \cK_i(U)=\Omega^\infty\Map^{\mathrm{Sp}}_{\SH(U)}(\one_U,A_{i,U}),
\]
where $\Map^{\mathrm{Sp}}$ denotes the mapping spectrum. They are
$\mathbb A^1$-invariant cdh sheaves by absoluteness and
\cite[Proposition~3.30]{BEM}, and $\cK_0\simeq\cK$ by effective-cover
adjunction. The first two slice triangles give fiber sequences
\[
 \cK_1\longrightarrow\cK_0\xrightarrow{\rk}\underline{\Z}_{\rm cdh},
 \qquad
 \cK_2\longrightarrow\cK_1\xrightarrow{\det}B_{\rm cdh}\Gm .
\]
These identifications follow from the rank and first Chern class
comparisons \cite[(4.41)--(4.42) and Corollary~9.12(1)]{BEM}.
Here $B_{\rm cdh}\Gm$ is the cdh sheafification of the space of line
bundles. The sections are $k\mapsto k[\mathcal O]$ and
$L'\mapsto[L']-[\mathcal O]$, the latter with its canonical
identification $\det([L']-[\mathcal O])\simeq L'$.
For either projection $u$ with pointed section $v$, subtraction
$y\mapsto y-v(u(y))$ lands in the homotopy fiber and restricts
there to the identity, since $u$ is additive. Hence $\cK_b$ is a
pointed retract of $\cK_0$ for $b=0,1,2$. Rationalization on its
diagonal cohomology is injective by the universal calculation in
Lemma~\ref{lem:arithmetic-evaluation}.

Use the window of length $n-b$. On $\cK_b$ consider the operations
\[
 \frac{L}{T_{n-b}}F^{(n-b)}_{b,b},
 \qquad
 \frac{L}{n!}\,s_n\circ (\cK_b\to\cK_0).
\]
Their common rationalization is $L\ch_n$ on $\cK_b$, so injectivity
identifies them integrally. On $Y$, evaluate them on
the image in $A_{b,Y}$ of the lift $\hat z$ to $A_{c,Y}$, taking
$\pi_s$ for $s>0$. On $X$, apply Lemma~\ref{lem:arithmetic-evaluation}
with $r=c$, $a=b$ and $M=n-b$, hence $j=n-c$. The square in the
proof of Theorem~\ref{thm:arithmetic-rr}, with these indices, gives
\[
 \frac{L}{n!}\chi_{n,s}(f_!x)
 =f_*\!\left(L[\ch(x)\Td(T_f)]^{(s)}_{n-c}\right).
\]
For $n<c$ the upper horizontal map is zero by slice orthogonality.
\end{proof}

For $n=c=7$, this argument gives $L=10080=2\cdot7!$. Thus it gives
a factor of $2$ on the Newton class for a general target. On smooth
arithmetic schemes, cohomology with supports recovers the coefficient
$1$ in every source-weight-zero case.

\begin{corollary}\label{cor:arithmetic-smooth}
Let $f\colon X\to Y$ be a projective morphism between smooth
quasi-projective $S$-schemes of pure relative dimensions. Then the
formulas of Theorem~\ref{thm:main} hold in $H^{2n-s}(Y,\Z(n))$
for $x\in K_s(X)$, with $d=\rk T_f$ and motivic Gysin map $f_*$. They have multiplier
$C_{n,d}$ for $m\ge0$ and give $\chi_{n,s}(f_!x)=0$ for $m<0$.
Here $f_!$ is induced by derived direct image of perfect complexes.
\end{corollary}
\begin{proof}
The schemes are regular, and a projective embedding into a projective
space over $Y$ gives a regular immersion. The identification of the
$\KGL$-Gysin map with the usual direct image in this case is
\cite[Example~3.3.4(3)]{DegliseOrientation}, together with
\cite[\S4.4.3]{DJK}. For $m>0$ the assertion follows from
Theorem~\ref{thm:arithmetic-rr} if $d\ge0$ and from
Corollary~\ref{cor:arithmetic-shift} if $d<0$.

For the remaining cases let $Z=f(X)_{\rm red}$. For a pure
codimension-$c$ closed subset $D\subset Y$, localization identifies
the cohomology with supports in bidegree $(2n-s,n)$ with degree $-s$
hypercohomology of the Zariski sheafified cycle complex
$z^{n-c}(D,\bullet)$, placed in cohomological degrees $-\bullet$
\cite[\S4.3.1, (4.6)]{FramedTransfers}. If $c>n$ this complex is
zero. For closed subsets all of whose components have codimension
$>n$ the same follows by stratification and localization for nested
supports. In particular, restriction off any such subset is injective.
Since $f_!x$ restricts to zero off $Z$ and
$\operatorname{codim}_Y Z\ge -d$, this proves the assertion for $m<0$.

Suppose $m=0$ and $n>0$. By this injectivity, discard from $Y$ the
image components of codimension $>n$ and the intersections of the
remaining components. For a remaining component the supported
complex is $z^0(Z,\bullet)$. In each simplicial degree its
generators are the irreducible components times the simplex, and
every face induces the identity. Its normalized sheaf complex is
concentrated in degree zero. Thus $H_Z^{2n-s}(Y,\Z(n))=0$ for $s>0$.
Localization gives $\chi_{n,s}(f_!x)=0$ in this case.

For $s=0$, each remaining component of $Z$ has the dimension
of $X$ and is horizontal since $X$ is smooth over $S$.
Over its generic point the map is finite with separable residue
extensions, since the source is reduced and the generic fields have
characteristic zero. Shrinking about these generic points removes
another subset of codimension $>n$ in $Y$. We can therefore assume
that $Z$ is regular and $f=i\circ g$, with $g\colon X\to Z$ finite
\'etale and $i\colon Z\hookrightarrow Y$ a regular immersion of
codimension $n$. Here properness allows us to remove the image of
the non-quasi-finite or non-\'etale locus. For $y=g_!x$,
Riemann--Roch without denominators for $i$ gives $c_j(i_!y)=0$ for
$0<j<n$ and $c_n(i_!y)=(-1)^{n-1}(n-1)!\,i_*\rk(y)$
\cite[Theorem~3.1]{NavarroHigherNoDen}. Newton's identity and
\eqref{eq:arithmetic-rr} in weight zero for $g$ now give
\[
 \chi_{n,0}(f_!x)=n!\,i_*\rk(g_!x)=n!\,f_*\rk(x).
\]
Injectivity of restriction gives the identity on $Y$. If $m=n=0$,
Theorem~\ref{thm:arithmetic-rr} already applies with $T_0=1$.
\end{proof}

For the same smooth arithmetic schemes, put $m=n+d\ge0$ and
$M=\max(n,m)$. The absolute Todd formula reads
\[
 T_M[\ch(f_!x)\Td(T_{Y/S})]^{(s)}_n
 =f_*\!\left(T_M[\ch(x)\Td(T_{X/S})]^{(s)}_m\right).
\]
Indeed, put $M_i=\max(i,i+d)$. Since $M_i+n-i=M$,
$T_{M_i}T_{n-i}\mid T_M$ by \eqref{eq:divisibility}. Multiply
\eqref{eq:arithmetic-rr} in weight $i$ by
$T_M D_{n-i}(T_{Y/S})/(T_{M_i}T_{n-i})$ and sum over $0\le i\le n$.
If $i+d<0$, the target term vanishes by
Corollary~\ref{cor:arithmetic-smooth}. If $d>0$, the missing
source terms have index $0\le j<d$. Each is the integral product
\[
 \frac{T_M}{T_jT_{m-j}}\,
 \bigl(T_j[\ch(x)\Td(T_f)]^{(s)}_j\bigr)\,
 f^*D_{m-j}(T_{Y/S}).
\]
The parenthesized class pushes forward to weight $j-d<0$, so the
projection formula gives zero. The full convolution now agrees with
the asserted integral polynomial by
$\Td(T_f)f^*\Td(T_{Y/S})=\Td(T_{X/S})$.

\section{Examples and optimal coefficients}\label{sec:examples}

We first write out the integral Chern character in weights zero, one and two.
We then recover the classical Riemann--Roch formula for a curve and check
Theorem~\ref{thm:main} on twists of pullback classes on the relative
projective line. Finally, we use finite covers to examine which
coefficients are necessary.

The first Todd numbers and top Newton coefficients are
\[
\begin{array}{c|rrrrr}
 q&0&1&2&3&4\\\hline
 T_q&1&2&12&24&720\\
 T_q/q!&1&2&6&4&30.
\end{array}
\]
Recall from Theorem~\ref{thm:intro-stable} that the integral map
$\Phi_M\colon W_M\to H_M$ induces on $K_s(X)$ the operation whose
weight-$q$ component is $(T_M/q!)\chi_{q,s}$. For the truncation
$W_2$, we have $T_2=12$, so on $K_0(X)$ this gives
\[
 x\longmapsto
 \bigl(12\rk(x),\ 12c_1(x),\ 6(c_1(x)^2-2c_2(x))\bigr).
\]
On $K_s(X)$, $s>0$, it is
$x\mapsto(0,12c_{1,s}(x),-12c_{2,s}(x))$.
These formulas illustrate the common multiplier $T_M$ of a finite
truncation and the linearization on positive simplicial suspensions.

For a smooth projective geometrically connected curve $C/k$ of
genus $g$ and a vector bundle $V$ on $C$, the case $n=0$, $d=1$ gives
\[
 2\chi(C,V)=2\deg V+(2-2g)\rk(V).
\]
Here the target is $\Z$, so division by $2$ recovers
$\chi(C,V)=\deg V+(1-g)\rk(V)$.

For the projection $p\colon\mathbf P^1_Y\to Y$, take $a\in K_s(Y)$ and set
$x=p^*a\,[\mathcal O(\ell)]$ for $\ell\in\Z$.
Writing $h=c_1(\mathcal O(1))$, one has
\[
 h^2=0,\qquad\Td(T_p)=1+h,\qquad p_!x=(\ell+1)a.
\]
The integral Newton identity is
\[
 \chi_{j,s}(x)=p^*\chi_{j,s}(a)
       +j\ell h\,p^*\chi_{j-1,s}(a)\qquad(j\ge1).
\]
For $s=0$ this follows from the Chern-root formula for tensoring
with a line bundle; pulling back the universal identity to a
simplicial suspension gives the same formula for $s>0$.
Using \eqref{eq:integral-polynomial}, both sides of
Theorem~\ref{thm:main}, with $C=T_{n+1}$, are consequently
\[
 (\ell+1)\frac{T_{n+1}}{n!}\chi_{n,s}(a).
\]
This calculation takes place integrally: the two coefficients on
the pushforward side are $T_{n+1}/n!$ and
$\ell T_{n+1}/n!$. In particular it applies to all higher
$K$-groups, with $\chi_{0,s}=0$ for $s>0$.

\subsection{Prime obstructions and finite covers}

\begin{proposition}
\label{prop:torsor-prime-obstructions}
Let $n\ge1$, let $p$ be a prime invertible in $k$, and choose $N\ge n$.
There is a finite \'etale $\mu_p$-torsor
$f_p:X\to Y_p$ between smooth quasi-projective $k$-schemes such that
\[
 \CH^n(Y_p)=\Z/p\cdot h^n,
 \qquad
 \chi_{n,0}\bigl((f_p)_![\mathcal O_X]\bigr)
 =\begin{cases}
   -h^n,&p-1\mid n,\\
   0,&p-1\nmid n.
  \end{cases}
\]
Consequently, in characteristic zero, every integer $a$ for which
\[
 a\chi_{n,0}(f_!x)=f_*\bigl(a\chi_{n,0}(x)\bigr)
\]
holds for all finite \'etale $f$ and all $x\in K_0(X)$ is divisible by
\[
 b_n^{\mathrm{cyc}}=\prod_{p-1\mid n}p.
\]
In particular, $C_{n,0}=T_n$ is optimal in this Newton normalization
for $n=1,2,4$.
\end{proposition}

\begin{proof}
Put
\[
 Y_p=\mathcal O_{\mathbf P^N}(-p)^\times,
 \qquad
 X=\mathcal O_{\mathbf P^N}(-1)^\times
       \simeq\mathbf A^{N+1}-\{0\},
\]
and define $f_p$ by taking the $p$th tensor power of a nonzero vector
in the tautological line. This is a $\mu_p$-torsor, hence finite
\'etale when $p$ is invertible in $k$; it requires no roots of unity
in $k$. Let $L$ be the pullback of $\mathcal O_{\mathbf P^N}(1)$
to $Y_p$, and put $h=c_1(L)$. Localization at the zero section gives
\[
 \CH^*(Y_p)=\Z[h]/(h^{N+1},ph),
 \qquad
 (f_p)_*\mathcal O_X=\bigoplus_{j=0}^{p-1}L^j.
\]
The second identity is the natural grading of the algebra of the
torsor, using the tautological trivialization of $L^p$.
The integral Newton polynomial therefore gives
\[
 \chi_{n,0}\bigl((f_p)_![\mathcal O_X]\bigr)
 =\left(\sum_{j=0}^{p-1}j^n\right)h^n.
\]
For $n\ge1$, the power sum is $-1$ modulo $p$ if $p-1\mid n$
and is zero otherwise, by summing over the cyclic group
$\mathbf F_p^\times$. Since $h^n$ has order $p$, this proves the
first assertion. Now $T_{f_p}=0$ and
$\chi_{n,0}([\mathcal O_X])=0$, so every prime occurring in $b_n^{\mathrm{cyc}}$
must divide $a$.

Finally,
\[
 (b_1^{\mathrm{cyc}},b_2^{\mathrm{cyc}},b_4^{\mathrm{cyc}})=(2,6,30)
 =\left(\frac{T_1}{1!},\frac{T_2}{2!},\frac{T_4}{4!}\right).
\]
The upper bound of Theorem~\ref{thm:main} supplies these coefficients,
and the necessary divisibilities prove their optimality.
\end{proof}

For $n=3$ this family forces only $2\mid a$, whereas the stable
construction gives $T_3/3!=4$. These cyclic examples detect $b_n^{\mathrm{cyc}}$; in general this
is smaller than $L_n=\operatorname{rad}(T_n/n!)$, as already occurs
for $n=5$, where the two integers are $2$ and $6$.

\begin{proposition}
\label{prop:finite-etale-optimal}
Let $n\ge1$ and $L_n=\operatorname{rad}(T_n/n!)$.
For every finite \'etale morphism $f\colon X\to Y$ of smooth
quasi-projective varieties in characteristic zero and every
$x\in K_0(X)$,
\[
 L_n\bigl(s_n(f_!x)-f_*s_n(x)\bigr)=0
 \quad\text{in }\CH^n(Y).
\]
No smaller positive coefficient works uniformly, already over
$\mathbf C$.
\end{proposition}

\begin{proof}
The identity is the covering theorem of Fulton--MacPherson,
in the form recorded in \cite[\S2.d.1, equation~(2.7)]{Pappas}.
For completeness we describe how the lower bound passes from
group cohomology to Chow groups.
By \cite[Theorem~III]{EvensKahn}, a coefficient valid for all
induced representations of finite groups must be divisible by
$L_n$. Let $H\subset G$ and an $H$-representation $V$ detect
one of the prime divisors. Take sufficiently many copies of the
regular complex representation of $G$, and let $U$ be their free
locus. The complement has complex codimension greater than $n$.
Consequently $U(\mathbf C)$ is $2n$-connected, and the classifying
map $U/G\to BG$ induces an isomorphism in integral cohomology
through degree $2n$. The same holds for $H$.

The quotient map $f\colon U/H\to U/G$ is finite \'etale between
smooth quasi-projective varieties. The associated bundle
$E=U\times^H V$ has direct image associated with
$\operatorname{Ind}_H^G V$. Integral cycle classes commute with
Chern classes and with finite proper pushforward. Thus the Chow
class $s_n(f_*E)-f_*s_n(E)$ maps to the corresponding nonzero
group-cohomological defect. Any universal annihilator of these
Chow classes must therefore contain every prime dividing $L_n$.
This proves minimality.
\end{proof}

\begin{remark}
In the multiplier convention $C/n!$ of Theorem~\ref{thm:main},
the optimal finite-\'etale multiplier on $K_0$ is $n!L_n$.
This equals $T_n$ for $n=1,2,4$, but is $12$ rather than $24$
for $n=3$. The result does not determine the optimum for arbitrary
projective morphisms of relative dimension zero, whose relative
tangent class need not vanish.
\end{remark}

\subsection{Higher degree}

\begin{remark}
The optimal multiplier for the stable operation is a different
question from the optimal multiplier for a particular direct-image
formula on $K_s$. For example, if $f\colon X\to Y$ is finite locally
free, then under $\CH^1(-,1)=\Gamma(-,\mathcal O^\times)$ one has
\[
 c_{1,1}(f_!x)=\operatorname{Nm}_{X/Y}(c_{1,1}(x)).
\]
Indeed, $c_{1,1}$ is the determinant and restriction of scalars
satisfies $\det_Y(f_*g)=\operatorname{Nm}_{X/Y}(\det_Xg)$ for an
automorphism $g$ of a vector bundle; this is the determinant-norm
identity for a finite locally free algebra, and it passes to $K_1$.
The formula of Theorem~\ref{thm:main} supplies its multiple by $2$.
More generally, $\chi_{n,s}=(-1)^{n-1}n c_{n,s}$ for $s>0$,
so a bound for the stable Newton operation need not be sharp for
the individual higher Chern class. The present argument gives the
displayed universal bounds; the torsion examples above concern $K_0$.
\end{remark}

\begin{example}
Let $Y_2=\mathcal O_{\mathbf P^2}(-2)^\times$, and take the double
cover $f_2:\mathbf A^3-0\to Y_2$ of
Proposition~\ref{prop:torsor-prime-obstructions}, over a field of
characteristic different from $2$. Write $L=\pi^*\mathcal O(1)$
and $h=c_1(L)$. Base change the cover to
\[
 g:(\mathbf A^3-0)\times\mathbf G_m\longrightarrow
       Y_2\times\mathbf G_m.
\]
Let $\eta=[t]\in K_1(\mathbf G_m)$ denote the coordinate unit and
its pullback to the target, and put $x=g^*\eta$.
The projection formula gives $g_!x=(1+[L])\eta$.

For a line bundle $M$, put $v=c_1(M)$. Tensoring a virtual class
$z$ of rank zero satisfies the integral identity
\[
 c_t([M]z)=c_{\,t/(1+vt)}(z).
\]
It follows by expressing $z$ as the difference of two bundles of
equal rank and comparing their Chern roots. In degree $n$ this is
an identity of integral universal Chern polynomials; pulling it
back to a positive simplicial suspension yields
\begin{equation}\label{eq:higher-chern-line-tensor}
 c_{n,s}([M]z)=\sum_{j=1}^n
 (-1)^{n-j}\binom{n-1}{j-1}v^{n-j}c_{j,s}(z)
 \qquad(s>0).
\end{equation}
The higher Chern classes are additive, since the cross terms in
the Whitney sum formula vanish on a positive suspension.
Moreover,
\[
 \CH^j(\mathbf G_m,1)
 \simeq\CH^j(k,1)\oplus\CH^{j-1}(k)=0\quad(j>1),
 \qquad c_{1,1}(\eta)=\{t\}.
\]
Formula~\eqref{eq:higher-chern-line-tensor} and additivity now give
directly
\[
 c_{3,1}(g_!x)=h^2\{t\},
 \qquad \chi_{3,1}(g_!x)=3h^2\{t\}.
\]
The first class is killed by $2$, and the localization boundary
at $t=0$ maps it to the nonzero class $h^2\in\CH^2(Y_2)$.
Thus both displayed classes are nonzero of order $2$.
Since $T_g=0$ and $c_{3,1}(x)=\chi_{3,1}(x)=0$, the unmultiplied
direct-image identity fails for the individual higher Chern class
as well as for the Newton class. In Newton normalization any
universal coefficient must be even, while Theorem~\ref{thm:main}
supplies $T_3/3!=4$.
\end{example}

\begin{question}[Optimal coefficients in higher degree]
\label{question:higher-optimal}
In characteristic zero, for fixed $n\ge1$ and $s>0$, determine the
least positive integer $a_{n,s}$ for which
\[
 a_{n,s}\bigl(\chi_{n,s}(f_!x)-f_*\chi_{n,s}(x)\bigr)=0
\]
for every finite \'etale morphism between smooth quasi-projective
schemes and every $x\in K_s(X)$. Theorem~\ref{thm:main} gives
$a_{n,s}\mid T_n/n!$. The determinant-norm identity gives
$a_{1,1}=1$, and the preceding example gives
$2\mid a_{3,1}\mid4$. Is $a_{3,1}=2$?
\end{question}

The comparison with the classical covering theorem in positive
$K$-degrees is not established here.
The analogous question for $c_{n,s}$ has the a priori upper bound
$T_n/(n-1)!$. The universal line-tensor calculation above detects
a nonzero higher Chern defect without cancelling the factor $n$
in $\chi_{n,s}=(-1)^{n-1}n c_{n,s}$.

\section*{Funding}
Both authors were supported by FEDER/AEI/MICINN under Grant No.~PID2022-142024NB-I00 (SIGA).

\section*{Use of AI tools}
We used ChatGPT (OpenAI) and Claude (Anthropic) to adapt the Bockstein
calculation of Dundas and \O{}stv{\ae}r~\cite{DO} in
Section~\ref{sec:integral-realization}. We have checked the resulting
computations. 

\needspace{10\baselineskip}

\par\bigskip
\needspace{10\baselineskip}
\begingroup
\small
\setstretch{1}
\noindent\paperemails\par
\endgroup
\end{document}